\documentclass[12pt]{article}
\usepackage{amsmath}
\usepackage{amsfonts}
\usepackage{amssymb}
\usepackage{latexsym}
\usepackage{graphicx}
\usepackage{enumerate}
\usepackage[shortlabels]{enumitem}
\usepackage{mathtools}
\usepackage{hyperref}

\usepackage{amsmath, amssymb,latexsym}
\usepackage{color}
\usepackage{xcolor}
\usepackage{fullpage}

\def\R{\mathbb{R}}
\def\C{\mathbb{C}}
\def\F{\mathbb{F}}
\def\la{\lambda}

\def\Z{\mathbb{Z}}
\def\K{\mathbb{K}}

\DeclareMathOperator{\rank}{rank}

\newtheorem{theorem}{Theorem}[section]
\newtheorem{proposition}[theorem]{Proposition}

\newtheorem{lemma}[theorem]{Lemma}
\newtheorem{definition}[theorem]{Definition}
\newtheorem{corollary}[theorem]{Corollary}
\newtheorem{remark}[theorem]{{\sc Remark}}

\title{Revisiting the matrix polynomial greatest common divisor}
\author{
Vanni Noferini\thanks{Aalto University, Department of Mathematics and Systems Analysis, P.O. Box 11100, FI-00076, Aalto, Finland. (\texttt{vanni.noferini@aalto.fi}). Supported by an Academy of Finland grant (Suomen Akatemian p\"{a}\"{a}t\"{o}s 331240).}
\and 
Paul Van Dooren\thanks{Universit\'{e} catholique de Louvain, Department of Mathematical Engineering, Av. Lemaitre 4, B-1348
Louvain-la-Neuve, Belgium. (\texttt{vandooren.p@gmail.com}). Supported by an Aalto Science Institute Visitor Programme.}
}
\begin{document}
\maketitle
\begin{abstract}
        In this paper we revisit the greatest common right divisor (GCRD) extraction from a set of  polynomial matrices $P_i(\lambda)\in \F[\la]^{m_i\times n}$,  $i=1,\ldots,k$ with coefficients in a generic field $\F$, and with common column dimension $n$. We give necessary and sufficient conditions for a matrix $G(\la)\in \F[\la]^{\ell\times n}$ to be a GCRD using the Smith normal form of the $m \times n$ compound matrix $P(\lambda)$ obtained by concatenating $P_i(\lambda)$ vertically,
where $m=\sum_{i=1}^k m_i$. We also describe the complete set of degrees of freedom for the solution $G(\la)$, and we link it to the Smith form and Hermite form of $P(\la)$. We then give an algorithm for constructing a particular minimum rank solution for this problem when $\F=\C$ or $\R$, using state-space techniques. This new method works directly on the coefficient matrices of $P(\la)$, using orthogonal transformations only.
The method is based on the staircase algorithm, applied to a particular pencil derived from a generalized state-space model of $P(\la)$. 
\end{abstract}
\\
\\
\textbf{AMS subject classifcation.} 15A22, 15A24, 65F45, 93-08, 93B52

\section{Introduction} \label{sec:intro}

The notions of a common right divisor (CRD) and of a greatest common right divisor (GCRD) of a set of polynomial matrices is a natural extension of the analogous concepts for scalar polynomials. As commutativity is lost in the matrix case, one needs first to specify on which side (right or left) the common divisor should lie. The following definition of CRDs and GCRDs is classical \cite{GLR82,kailath,wolovich}; while \cite{GLR82,kailath} only consider the field of complex numbers $\C$, the same definition is in fact still valid when $\F$ is any arbitrary field \cite{wolovich}.

\begin{definition}[Common right divisor and greatest common right divisor] \label{def:GCRD}
Let $\{P_i(\la)\}_{i=1}^k$, with $P_i(\la)\in \F[\la]^{m_i\times n}, \;i=1,\ldots,k$, be a non-empty set of polynomial matrices all having the same number of columns. We say that $D(\la) \in \F[\la]^{p \times n}$ is a common right divisor (CRD) of the set $\{P_i(\la)\}_{i=1}^k$ if it is a right divisor of each $P_i(\la)$, i.e., if there exist polynomial matrices $Q_i(\la)\in \F[\la]^{m_i\times p}, \; i=1,\ldots,k$ such that $P_i(\la)=Q_i(\la) D(\la)$ for all $i=1,\ldots,k$. Furthermore, $G(\la) \in \F[\la]^{\ell \times n}$ is said to be a greatest common right divisor (GCRD) of the set $\{P_i(\la)\}_{i=1}^k$ if, for every CRD $D(\la) \in \F[\la]^{p \times n}$ of the same set, $D(\la)$ is a right divisor of $G(\la)$, i.e., there exist a polynomial matrix $Q(\la) \in \F[\la]^{\ell \times p}$ such that $G(\la)=Q(\la)D(\la)$. 
\end{definition}

Even though Definition \ref{def:GCRD} applies to any non-empty set of polynomial matrices $P_i(\la)\in \F[\la]^{m_i\times n}$ with the same number of columns $n$, the authors of \cite{GLR82,kailath,wolovich} immediately restrict their discussion of the problem by assuming (at least) that the compound matrix 
\begin{equation} \label{eq:compound} 
P(\la):=
 \left[\begin{array}{c} P_1(\la) \\ \vdots \\ P_k(\la) \end{array}\right]\in \F[\la]^{m\times n}, \quad m:=\sum_{i=1}^k m_i,
\end{equation}
has full normal rank $n$ and the GCRD is square (and therefore $n\times n$) and nonsingular. This simplification is crucial in the design of several of the existing algorithms, since they do not extend to the rank deficient case. Once a GCRD has been found, the polynomial factors 
$Q_i(\la)$ of Definition \ref{def:GCRD} are said to be {\em right coprime} \cite{amz}, meaning that they only have common right divisors with a trivial Smith form (i.e., with elementary divisors 1). In the case where \eqref{eq:compound} has full normal rank, this in turn implies that the square CRDs of the coprime matrices $Q_i(\la)$ must be unimodular.

As mentioned above, there is obviously a dual problem of finding a greatest common left divisor (GCLD) for a set of polynomial matrices all having the same number of rows, but since it is obvious that the latter can be reduced to the GCRD extraction simply by transposing all the defining equations, we restrict ourselves in this paper to the right factor extraction. 
One of the main applications of GCRD and GCLD extraction (for $\F=\C$ or $\R$) is the reduction of polynomial system quadruples $\{A(\la),B(\la),C(\la),D(\la)\}$, introduced by Rosenbrock
\cite{Rosenbrock} in linear system theory, and recently become of interest also in the context of nonlinear eigenvalue problems, see e.g. \cite{dmqv,dqv,nnpq} and the references therein. Such system quadruples are a minimal representation of the (rational)
transfer function $R(\la):=C(\la)A(\la)^{-1}B(\la)+D(\la)$
provided the pairs $\{-A(\la),C(\la)\}$ and $\{-A(\la),B(\la)\}$
are right and left coprime, respectively. In these 
system quadruples, $A(\la)$ is always assumed to be invertible over the field of fractions $\F(\la)$, implying that the compound matrices
$ \left[\begin{array}{cc}-A(\la)\\ C(\la)\end{array}\right],  \left[\begin{array}{cc}-A(\la)& B(\la)\end{array}\right]
$
have automatically full column and row rank, respectively. If the system quadruple $\{A(\la),B(\la),C(\la),D(\la)\}$ is not minimal, one can obtain a minimal one by extracting a GCRD of the pair 
$\{-A(\la),C(\la)\}$, followed by the computation of a GCLD of the (updated) pair $\{-A(\la),B(\la)\}$ \cite{coppel}. Other relevant applications, for which the compound matrix \eqref{eq:compound} may also be rank deficient, include nonsingular factorizations of polynomial matrices \cite{emre, moness} or coprimeness tests for polynomial matrices \cite{moness}.
It is therefore not surprising that the GCRD and GCLD extraction problem
has been extensively studied in the systems theory literature
\cite{kailath}. Furthermore, the problem of computing a GCRD (or GCLD) is also mathematically interesting per se and has been well studied in the matrix theory literature, with particular emphasis on its connection with Sylvester or B\'{e}zout resultant matrices \cite{gklr81,gklr82,GLR82,lt,nnt}. Nevertheless, all the numerical algorithms that we could find in the literature are designed for the full rank case only \cite{Bitmead,emre,GLR82,kailath,SilV81}
and their extension to the rank deficient case is not discussed. We also mention that recently the related (but distinct) problem of computing an \emph{approximate} GCRD has received attention \cite{fgm}, but again under the full rank assumption.

The present manuscript has two goals. On one hand, we review and extend the theory of GCRDs, paying particular attention to the non-full rank case that apparently has received little to no attention in the existing literature. Our description is built on the Hermite and Smith normal forms, and thus is similar in spirit (while making much weaker assumptions) to the approach of the control theory literature \cite{kailath,wolovich} but significantly differs from the approach based on common restrictions described in \cite{GLR82} and the references therein. On the other hand, for $\F=\C$ or $\R$, we propose a new algorithm for the computation of a GCRD, that is based on generalized state-space realizations. The new algorithm is inspired by previous approaches based on state-space realizations \cite{emre,SilV81}, but it differs from them in (a) using generalized state-space realizations, thus being able to accept as input any non-empty set of polynomial matrices with the same number of columns with no further restriction and (b) relying on the staircase algorithm \cite{vd79} which is numerically stable.

The organization of the paper is as follows. In Section \ref{sec:theory},
we recall the theory of polynomial matrices with coefficients over any arbitrary field, needed to define correctly greatest common divisors and their connection to the Smith normal form and the Hermite normal form of $P(\la)$. We also develop a comprehensive list of properties that such greatest common right divisors must have and how any generic GCRD can be derived from a particular, minimal-size, solution that we label a \emph{compact GCRD}. We then restrict our attention to the case $\F=\C$ or $\R$, with the aim of developing a numerical algorithm. In Section \ref{sec:feedback} we recall useful properties of state-space realizations and feedback and show that feedback can be interpreted as a factorization of a given transfer function. In Section \ref{sec:gssm} we then extend these ideas to generalized state space realizations of a general polynomial matrix $P(\la)$ and show how to use this to compute a compact GCRD of $P(\la)$. Section \ref{sec:numerics} then discusses numerical aspects of GCRD
algorithms and gives numerical experiments illustrating our new algorithm.
We end this paper with some concluding remarks in Section \ref{sec:conclusions}.

\section{Theory of greatest common right divisors} \label{sec:theory}

In this section we review and extend the theory of GCRDs; some of the results that we give are to our knowledge new, while others appeared in the literature (see e.g. \cite{GLR82,kailath,wolovich} and the references therein) but sometimes under additional assumptions that we have relaxed or removed.
In particular, we establish the existence of a GCRD of any non-empty set of polynomial matrices with the same number of columns, as well as its uniqueness up to fixing the number of rows and up to left multiplication by a unimodular matrix. We start with an obvious, but useful, lemma. It relies on Definition \ref{def:GCRD}, given earlier.

\begin{lemma}\label{lem:conditionforg}
Let $\{P_i(\la)\}_{i=1}^k$, with $P_i(\la)\in \F[\la]^{m_i\times n}, \;i=1,\ldots,k$,  be a non-empty set of polynomial matrices and suppose that $G(\la) \in \F[\la]^{\ell \times n}$ is a GCRD of $\{ P_i(\la)\}_{i=1}^k$. Suppose moreover that $D(\la) \in \F[\la]^{p \times n}$ is a CRD of $\{ P_i(\la)\}_{i=1}^k$. If $G(\la)$ is a right divisor of $D(\la)$, then $D(\la)$ is a GCRD of $\{ P_i(\la)\}_{i=1}^k$.
\end{lemma}
\begin{proof}
By assumption there is a polynomial matrix $A(\la)$ such that $D(\la)=A(\la)G(\la)$. Let $C(\la)$ be any CRD of $\{ P_i(\la)\}_{i=1}^k$; then, since $G(\la)$ is a GCRD of the same set, $G(\la)=Q(\la)C(\la)$ for some polynomial matrix $Q(\la)$. It follows that $D(\la)=A(\la)Q(\la)C(\la)$, showing that $C(\la)$ divides $D(\la)$ and hence that $D(\la)$ is a GCRD of the given set.  \hfill
\end{proof}

The next result provides a lower bound on the row size of a right CRD.

\begin{lemma}\label{lem:rowsofcrd}
Let $\{P_i(\la)\}_{i=1}^k$, with $P_i(\la)\in \F[\la]^{m_i\times n}, \;i=1,\ldots,k$,  be a non-empty set of polynomial matrices and suppose that $D(\la) \in \F[\la]^{p \times n}$ is a CRD of $\{ P_i(\la)\}_{i=1}^k$. Then, $p \geq r$ where $r=\rank P(\la)$ is the normal rank of the compound matrix \eqref{eq:compound}.
\end{lemma}
\begin{proof}
Since $D(\la)$ is a CRD of the given set, there must exist a polynomial matrix $N(\la) \in \F[\la]^{m \times p}$ such that $P(\la)=N(\la)D(\la)$. It follows that $r = \rank P(\la) \leq \rank D(\la) \leq p$.  \hfill
\end{proof}

\noindent Lemma \ref{lem:rowsofcrd} justifies the following definition.
\begin{definition}
Let $\{P_i(\la)\}_{i=1}^k$, with $P_i(\la)\in \F[\la]^{m_i\times n}, \;i=1,\ldots,k$,  be a non-empty set of polynomial matrices and suppose that their compound matrix \eqref{eq:compound} has rank $r$. Then, any GCRD (resp. CRD) $G(\la) \in \F[\la]^{r \times n}$ is called a compact GCRD (resp. compact CRD) of the set $\{P_i(\la)\}_{i=1}^k$.
\end{definition}
From now on, we will assume that $r \geq 1$. This is essentially no loss of generality, for if $r=0$ then $P_i(\la)$ are all zero matrices, and the problem of extracting a GCRD is trivial: $G(\la)$ is a GCRD of theirs if and only if it is a zero matrix with the same number of columns as the $P_i(\la)$.
\begin{proposition}\label{prop:fullrowrank}
Let $\{P_i(\la)\}_{i=1}^k$, with $P_i(\la)\in \F[\la]^{m_i\times n}, \;i=1,\ldots,k$,  be a non-empty set of polynomial matrices and suppose that their compound matrix $P(\la)$ as in \eqref{eq:compound} has rank $r \geq 1$. If $G(\la) \in \F[\la]^{r \times n}$ is a compact GCRD of $\{P_i(\la)\}_{i=1}^k$, then $G(\la)$ has full row rank.
\end{proposition}
\begin{proof}
From the equation $P(\la)=N(\la)G(\la)$ we have $r \leq \rank G(\la)$. \hfill
\end{proof}

\noindent In a sense, it is natural to restrict to compact GCRDs because, as we will see below, any other GCRD can be constructed from them. However, it is conceivable that there are scenarios where the number $\ell$ of rows of a GCRD $G(\la)\in\F[\la]^{\ell \times n}$ is strictly greater than $r$; for instance, it may be natural to take $\ell=n$, thus studying square GCRDs, even if $n>r$. For this reason, we will also study the case $\ell>r$.

\subsection{Background on normal forms of polynomial matrices}

Recall that a square polynomial matrix is said to be regular if it is invertible over the field of fractions $\F(\la)$, that is, if its determinant is a nonzero polynomial; and it is said to be a unimodular matrix if it is invertible over the ring $\F[\la]$, that is, if its determinant is a nonzero constant. Moreover, a polynomial matrix $N(\la)$, not necessarily square, is called left (resp. right) invertible if there is a polynomial matrix $L(\la)$ such that $L(\la)N(\la)=I$ (resp. $N(\la)L(\la)=I$). More precisely, the previous property defines left (resp. right) invertibility over the ring $\F[\la]$. Occasionally, we will also refer to the weaker property of being left (resp. right) invertible over the field of fractions $\F(\la)$, which means that there is a \emph{possibly rational} matrix $L(\la)$ such that $L(\la)N(\la)=I$ (resp. $N(\la)L(\la)=I$).

To further develop the theory of GCRDs, we will make extensive use of the Smith normal form and the Hermite normal form of a polynomial matrix. We recall them below.

\begin{theorem}[Smith normal form]
For any polynomial matrix $P(\la) \in \F[\la]^{m \times n}$ there exist unimodular matrices $U(\la) \in \F[\la]^{m \times m}$ and $V(\la) \in \F[\la]^{n \times n}$ such that $P(\la)=U(\la) S(\la) V(\la)$ and 
$$ S(\la):= \left[\begin{array}{ccc|c}
s_1(\la) & & & \\ & \ddots & & 0 \\ & & s_r(\la) \\ \hline & 0 & & 0
\end{array}\right]\in \F[\la]^{m \times n}$$
contains $r$ monic polynomials $s_i(\la)$ on diagonal, and zeros for the rest. The elements $s_i(\la)$ form a divisibility chain 
$$ s_1(\la) | s_2(\la) | \ldots | s_r(\la)
$$
and $r$ is the normal rank of $P(\la)$.
The matrix $S(\la)$ is called the Smith normal form of $P(\la)$ and is uniquely determined by $P(\la)$ \cite[Theorem 1.14.1 and Normalization 1.14.4]{Friedland}.
\end{theorem}

Occasionally in this paper we will use a {\em compact} Smith form as well. It is obtained by just deleting the last $n-r$ rows of 
$V(\la)$ and the last $m-r$ columns of $U(\la)$, as well as by keeping only the top-left $r \times r$ submatrix of $S(\la)$, thus yielding the decomposition
 $P(\la)=\hat U(\la)\hat S(\la)\hat V(\la)$,
where $\hat U(\la)\in \F[\la]^{m\times r}$ is left invertible,
$\hat S(\la)\in \F[\la]^{r\times r}$ is square, regular, and in Smith normal form, and $\hat V(\la)\in \F[\la]^{r\times n}$ is right invertible. When $r=n$, the right matrix $\hat V(\la)=V(\la)$ is still square and unimodular.

\begin{theorem}[Hermite normal form] \label{th:Hermite}
For any polynomial matrix $P(\la)\in \F[\la]^{m \times n}$ there exist a unimodular matrix $U(\la) \in \F[\la]^{m \times m}$ such that $P(\la)=U(\la)H(\la)$ and $H(\la) \in \F[\la]^{m \times n}$ satisfies the following properties:
\begin{itemize}
    \item It is in row echelon form and its zero rows (if any) are located below its nonzero rows (if any);
    \item The pivot (leftmost nonzero element) of each nonzero row is a monic polynomial and it is located strictly to the right of the pivots of the rows above;
    \item All elements below a pivot are zero and all the elements above a pivot are monic polynomials of degree strictly less than the pivot below;
    \item The number of its nonzero rows is $r$, the normal rank of $P(\la)$. 
\end{itemize} The matrix $H(\la)$ is called the Hermite normal form of $P(\la)$ and is uniquely determined by $P(\la)$ \cite[Theorem 1.12.7 and Corollary 1.12.11]{Friedland}.
\end{theorem}\label{eq:compactHermite}
When $m>r$, the bottom $m-r$ rows of $H(\la)$ are necessarily zero and one can then rewrite the decomposition $P(\la)=U(\la)H(\la)$ in the {\em compact} form $P(\la)=\hat U(\la) \hat H(\la)$,
 where $\hat U(\la)\in \F[\la]^{m\times r}$ is obtained by just deleting the rightmost $m-r$ columns of $U(\la)$ (and is therefore left invertible), and $\hat H(\la)\in \F[\la]^{r\times n}$ is obtained by deleting the $m-r$ bottom rows of $H(\la)$, which are zero. Moreover, it follows from Theorem \ref{th:Hermite} that, if $r=\rank P(\la)=n$, then $\hat H(\la)$ is square upper triangular and invertible.  

\subsection{Properties of greatest common right divisors}

We are now ready to further develop our analysis of the problem of finding a GCRD of a set of polynomial matrices. First, we show that there is no loss of generality in restricting to the case where $m:=\sum_{i=1}^k m_i \geq n$.

\begin{lemma} 
Let $\{P_i(\la)\}_{i=1}^k$, with $P_i(\la)\in \F[\la]^{m_i\times n}, \;i=1,\ldots,k$,  be a non-empty set of polynomial matrices such that $m:=\sum_{i=1}^k m_i <n$. Moreover, let $\delta=n-m$ and define $P_{k+1}(\la):=0_{\delta \times n}$. Then, $G(\la) \in \F[\la]^{\ell \times n}$ is a GCRD of $\{P_i(\la)\}_{i=1}^k$ if and only if it is a GCRD of $\{P_i(\la)\}_{i=1}^{k+1}$.
\end{lemma}

\begin{proof}
Suppose first that $G(\la)$ is a GCRD of $\{ P_i(\la)\}_{i=1}^{k+1}$. Then for $i=1,\ldots,k+1$ there are $Q_i(\la)$ such that $P_i(\la)=Q_i(\la)G(\la)$, implying that $G(\la)$ is a CRD of $\{ P_i(\la)\}_{i=1}^k$. If it was not a GCRD of the latter set, then such a set would admit another CRD $D(\la) \in \F[\la]^{p \times n}$ that is not a right divisor of $G(\la)$. In particular, $P_i(\la)=R_i(\la) D(\la)$ for all $i=1,\ldots,k$. But clearly $P_{k+1}(\la)=0_{\delta \times p} \cdot D(\la)$, and hence $D(\la)$ is a CRD of $\{ P_i(\la)\}_{i=1}^{k+1}$ contradicting that $G(\la)$ is a GCRD of $\{ P_i(\la)\}_{i=1}^{k+1}$.

Conversely suppose that $G(\la)$ is a GCRD of $\{ P_i(\la)\}_{i=1}^k$. Then for $i=1,\ldots,k$ there are $Q_i(\la)$ such that $P_i(\la)=Q_i(\la)G(\la)$, and obviously $P_{k+1}(\la)=0_{\delta \times \ell} \cdot G(\la)$, thus showing $G(\la)$ is a CRD of $\{ P_i(\la)\}_{i=1}^{k+1}$. If we suppose for a contradiction that $G(\la)$ is not a GCRD of the latter set, then we could find a CRD $D(\la)$ of $\{ P_i(\la)\}_{i=1}^{k+1}$ that is not a right divisor of $G(\la)$. But, following the argument above, then such $D(\la)$ would also be a CRD of $\{ P_i(\la)\}_{i=1}^k$. Hence, we easily reach the contradiction that $G(\la)$ cannot possibly be a GCRD of $\{ P_i(\la)\}_{i=1}^k$ either.   \hfill
\end{proof}

Next, we characterize the right null space of a GCRD.

\begin{lemma}\label{lem:kernel}
Let $\{P_i(\la)\}_{i=1}^k$, with $P_i(\la)\in \F[\la]^{m_i\times n}, \;i=1,\ldots,k$,  be a non-empty set of polynomial matrices such that $m:= \sum_{i=1}^k m_i \ge n$, and define their compound matrix $P(\la)$ as in \eqref{eq:compound}.
If $G(\la) \in \F[\la]^{\ell \times n}$ is a CRD of $\{P_i(\la)\}_{i=1}^k$ then $\ker G(\la)\subseteq \ker P(\la)$. If $G(\la) \in \F[\la]^{\ell \times n}$ is a GCRD of $\{P_i(\la)\}_{i=1}^k$ then $\ker G(\la)=\ker P(\la)$.
\end{lemma}

\begin{proof}
It is clear that $\ker G(\la) \subseteq \ker P(\la)$ whenever $G(\la)$ is a CRD of $\{P_i(\la)\}_{i=1}^k$. Indeed, due to the relation $P_i(\la)=Q_i(\la)G(\la)$ which is valid for $i=1,\ldots,k$, we have $v(\la) \in \ker G(\la) \Rightarrow v(\la) \in \ker P_i(\la) \ \forall \ i \Rightarrow v(\la) \in \ker P(\la)$. Now suppose that $G(\la)$ is a GCRD and consider a compact Smith decomposition $P(\la)=\hat U(\la)\hat S(\la)\hat V(\la)$ where $\hat U(\la)\in \F[\la]^{m\times r}$ is left invertible, $\hat S(\la)\in \F[\la]^{r \times r}$ is in Smith form, and $\hat V(\la) \in \F[\la]^{r \times n}$ is right invertible. By construction $\hat S(\la)\hat V(\la)$ is a CRD of $\{P_i(\la)\}_{i=1}^k$, and hence it is a right divisor of $G(\la)$. It follows that $\rank G(\la) \leq \rank \hat S(\la)\hat V(\la)=\rank \hat S(\la)=\rank P(\la)$. Hence, $\dim \ker G(\la)=\dim \ker P(\la)$ and therefore $\ker G(\la)=\ker P(\la)$.    \hfill
\end{proof}

We now state Theorem \ref{thm:gcdandhermite}, one of the main results of this section, implying among other things the existence of a GCRD of any given non-empty set of polynomial matrices all with the same number of columns (thus far in our analysis we have merely assumed that a GCRD always exists, taking the risk to have been pondering about the empty set). For convenience of exposition, however, we will first prove Theorem \ref{thm:gcdandhermite} under the assumption that \eqref{eq:compound} has full column rank, i.e., $n=r$ using the notation above. Later on, we will show how to remove this assumption.

\begin{theorem}\label{thm:gcdandhermite}
Let $\{P_i(\la)\}_{i=1}^k$, with $P_i(\la)\in \F[\la]^{m_i\times n}, \;i=1,\ldots,k$,  be a non-empty set of polynomial matrices such that $m:= \sum_{i=1}^k m_i \ge n$,  and let $P(\la)$ be their compound matrix as in \eqref{eq:compound}. Assume that $r=\rank P(\la) \geq 1$. Then, the following are equivalent:
\begin{enumerate}
    \item $G(\la) \in \F[\la]^{r \times n}$ is a compact GCRD of $\{P_i(\la)\}_{i=1}^k$;
        \item There exists a left invertible polynomial matrix $N(\la) \in \F[\la]^{m \times r}$ such that $P(\la)=N(\la)G(\la)$;
    \item There exists a unimodular matrix $U(\la) \in \F[\la]^{m \times m}$ such that 
    $$P(\la)=U(\la) \begin{bmatrix}
    G(\la)\\
    0
    \end{bmatrix}.$$
\end{enumerate}
Moreover, any compact GCRD $G(\la)$ of $\{P_i(\la)\}_{i=1}^k$ can always be expressed as a (left) linear combination $$G(\la)=\sum_{i=1}^k L_i(\la)P_i(\la)$$ for some polynomial matrices $L_i(\la) \in \F[\la]^{r \times m_i}$; equivalently, there is a polynomial matrix $L(\la) \in \F[\la]^{r \times m}$ such that $G(\la)=L(\la)P(\la).$ 
\end{theorem}

\begin{proof}[Proof of Theorem \ref{thm:gcdandhermite} in the case where $r=n$.]
\begin{itemize}
\item[$1 \Rightarrow 2$] Since $G(\la)$ is a CRD of $\{P_i(\la)\}_{i=1}^k$, there must exist matrices $N_i(\la)$ satisfying $P_i(\la)=N_i(\la)G(\la)$ for each $i$. Define $N(\la) := \begin{bmatrix}N_1(\la)\\ \vdots\\ N_k(\la)
\end{bmatrix}$, then by construction $P(\la)=N(\la)G(\la)$. Consider the compact Hermite decomposition $N(\la)=U(\la)H(\la)$ where $U(\la)$ is left invertible and $H(\la)$ is square upper triangular. It follows that $H(\la)G(\la)$ is a CRD of $\{P_i(\la)\}_{i=1}^k$. Hence, there exists some polynomial matrix $A(\la)$ such that $G(\la)=A(\la)H(\la)G(\la)$. However, by Lemma \ref{lem:kernel} and by the assumption that $P(\la)$ has full column rank, $\ker G(\la)=\{0\}$ and hence $A(\la)H(\la)=I$. It follows that $H(\la)=A(\la)^{-1}$ is unimodular and thus $N(\la)$ is left invertible.
\item[$2 \Rightarrow 3$] By, e.g., \cite[Theorem 3.3]{amz}, since $N(\la)$ is left invertible, it has a trivial Smith form and hence it is completable to a unimodular matrix $U(\la)$. The sought implication follows immediately.
\item[$3 \Rightarrow 1$] Partition
$$U(\la) = \begin{bmatrix}N_1(\la) & U_1(\la)\\
\vdots & \vdots\\ N_k(\la) & U_k(\la)
\end{bmatrix}, \qquad U(\la)^{-1} = \begin{bmatrix}L_1(\la) & \ldots & L_k(\la)\\
V_1(\la) & \ldots & V_k(\la)
\end{bmatrix}.$$ 
Then, since $P(\la)=N(\la)G(\la)$ implies $P_i(\la)=N_i(\la)G(\la)$, it is clear that $G(\la)$ is a CRD of the set $\{P_i(\la)\}_{i=1}^k$. On the other hand, $$\begin{bmatrix}G(\la)\\0\end{bmatrix}=U(\la)^{-1}P(\la) \Rightarrow G(\la) = \sum_{i=1}^k L_i(\la) P_i(\la).$$ Hence, it is clear that if $D(\la)$ is any CRD of $\{P_i(\la)\}_{i=1}^k$ then it must also be a right divisor of $G(\la)$. Therefore, $G(\la)$ is a GCRD of $\{P_i(\la)\}_{i=1}^k$.
\end{itemize}
Finally we note that the proof above, and in particular the argument to show that condition 3 implies condition 1, also implies the last statement. \hfill
\end{proof}

Since a unimodular matrix $U(\la)$ and a square matrix $G(\la)$ satisfying condition 3 can always be found -- for example via the compact Hermite form of $P(\la)$ or the Smith form of $P(\la)$ -- the proof above establishes the existence of a GCRD in the case where the compound matrix \eqref{eq:compound} has full column rank. We now proceed to explore the case where the compound matrix does not have full column rank $n$. In particular, finding compact GCRDs in this case can be reduced to the full column rank case.

\begin{theorem}\label{thm:kailathisconfusing}
Let $\{P_i(\la)\}_{i=1}^k$, with $P_i(\la)\in \F[\la]^{m_i\times n}, \;i=1,\ldots,k$,  be a non-empty set of polynomial matrices such that $m:= \sum_{i=1}^k m_i \ge n$,  and let $P(\la)$ be their compound matrix as in \eqref{eq:compound}. Suppose that $V(\la) \in \F[\la]^{n \times n}$ is a unimodular matrix such that $P(\la)V(\la)=\begin{bmatrix}Q(\la) & 0_{m \times (n-r)}
\end{bmatrix}$ where $$Q(\la)=\left[\begin{array}{c} Q_1(\la) \\ \vdots \\ Q_k(\la) \end{array}\right] \in \F[\la]^{m \times r}$$ with $Q_i(\la) \in \F[\la]^{m_i \times r}$, $i=1,\ldots,k$, and $r=\rank P(\la)=\rank Q(\la)$. Then, $D(\la) \in \F[\la]^{r \times r}$ is a compact GCRD of $\{Q_i(\la)\}_{i=1}^k$ if and only if
$G(\la)$ is a compact GCRD of $\{P_i(\la)\}_{i=1}^k$,
where
\begin{equation}\label{eq:gcdnotfullrank}
    G(\la)=U(\la)\begin{bmatrix}
D(\la)&0\\
\end{bmatrix}V(\la)^{-1} \in \F[\la]^{r \times n},
\end{equation} 
and $U(\la)\in \F[\la]^{r\times r}$ is an arbitrary unimodular polynomial matrix.
\end{theorem}
Before proving Theorem \ref{thm:kailathisconfusing}, we note that the existence of a unimodular matrix $V(\la)$ and a full column rank polynomial matrix $Q(\la)$ with the sought properties is clear, for example, from the Smith form of $P(\la)$.
\begin{proof}[Proof of Theorem \ref{thm:kailathisconfusing}.]
Suppose first that $D(\la)$ is a compact GCRD of $\{Q_i(\la)\}_{i=1}^k$ and write (by Theorem \ref{thm:gcdandhermite}) $Q(\la)=N(\la)D(\la)$ for some left invertible $N(\la) \in \F[\la]^{m \times r}$. Then, it is immediate that $P(\la)=N(\la)\begin{bmatrix}D(\la)&0\end{bmatrix}V(\la)^{-1}=N(\la)U(\la)^{-1}G(\la)$, so $G(\la)$ is a compact CRD of $\{P_i(\la)\}_{i=1}^k$. Now let $A(\la) \in \F[\la]^{p \times n}$ be any other CRD of the same set, so that $P(\la)=M(\la)A(\la)$ for some (not necessarily left invertible) $M(\la) \in \F[\la]^{m \times p}$. If we let $L(\la) \in \F[\la]^{r \times m}$ be a left inverse of $N(\la)$, then $D(\la)=L(\la)Q(\la)$ and therefore
 $$ G(\la) = U(\la) L(\la) \begin{bmatrix}Q(\la)&0_{m \times (n-r)}\end{bmatrix} V(\la)^{-1}=U(\la)L(\la)P(\la).$$ Hence, we have $G(\la)=U(\la)L(\la)M(\la)A(\la)$.
 This shows that $G(\la)$ is a GCRD of $\{P_i(\la)\}_{i=1}^k$.

For the reverse implication, we first need to show that any compact GCRD $G(\la)$ of $\{P_i(\la)\}_{i=1}^k$ must have the form \eqref{eq:gcdnotfullrank} for some polynomial matrix $D(\la)$ and some unimodular matrix $U(\la)$. By Lemma \ref{lem:kernel}, $\ker G(\la)=\mathrm{span}\left\{ V(\la) \begin{bmatrix} 0\\
I_{n-r}\end{bmatrix}\right\}$. Hence, necessarily such a compact GCRD $G(\la)$ must have the form
$$ G(\la) = \begin{bmatrix}
X(\la)&0\end{bmatrix}V(\la)^{-1} $$ for some appropriate polynomial matrix $X(\la) \in \F[\la]^{r \times r}$. However, obviously, we can write
$X(\la)=U(\la)D(\la)$ for some $r\times r$ unimodular matrix $U(\la)$ and some $r \times r$ polynomial matrix $D(\la)$ (just take $U(\la)$ arbitrarily and define $D(\la):=U(\la)^{-1}X(\la)$). Hence, we may assume that $G(\la)$ has indeed the form \eqref{eq:gcdnotfullrank}. 

On the other hand, since $G(\la)$ is a compact GCRD of $\{P_i(\la)\}_{i=1}^k$ then $P(\la)=M(\la)G(\la)$ for some $M(\la)\in \F[\la]^{n \times r}$. It follows that $D(\la)$ is a CRD of  $\{Q_i(\la)\}_{i=1}^k$, because \eqref{eq:gcdnotfullrank} then implies $Q(\la)=M(\la) U(\la) D(\la)$. Now let $A(\la) \in \F[\la]^{r \times r}$ be any compact GCRD of $\{Q_i(\la)\}_{i=1}^k$ (the existence of $A(\la)$ is guaranteed by the proof of Theorem \ref{thm:gcdandhermite} for the case of a full column rank compound matrix, since $Q(\la)$ has full column rank). By the first part of this proof
 $H(\la) := \begin{bmatrix} A(\la) & 0 \end{bmatrix} V(\la)^{-1} \in \F[\la]^{r \times n}$ is then a CRD of $\{P_i(\la)\}_{i=1}^k$. Hence, there is a polynomial matrix $W(\la) \in \F[\la]^{r \times r}$ such that $G(\la)=W(\la)H(\la)$. Therefore, $D(\la)=U(\la)^{-1}W(\la)A(\la)$,
which shows that $A(\la)$ is a right divisor of $D(\la)$, and hence $D(\la)$ is a GCRD of $\{Q_i(\la)\}_{i=1}^k$ by Lemma \ref{lem:conditionforg}. \hfill
\end{proof}

We are now ready to prove that Theorem \ref{thm:gcdandhermite} still holds even if we remove the assumption that the compound matrix $P(\la)$ has full column rank. By the comments that we have made above, establishing Theorem \ref{thm:gcdandhermite} in the most general case implies in particular that a compact GCRD exists for any given non-empty set of polynomial matrices all with the same number of columns.

\begin{proof}[Proof of Theorem \ref{thm:gcdandhermite} in the case where $\rank P(\la)=r<n$.]
In our proof for the case $r=n$, we only used such assumption in showing that condition 1 implies condition 2. Let us give another proof of the same implication, but this time assuming that $P(\la)$ has rank $r<n$. Then, by Theorem \ref{thm:kailathisconfusing}, $G(\la)$ has the form $U(\la) \begin{bmatrix}G_c(\la)&0\end{bmatrix}V(\la)^{-1}$ where $P(\la)V(\la)=\begin{bmatrix}Q(\la)&0\end{bmatrix}$, $Q(\la)$ has full column rank $r$ and, by Theorem \ref{thm:gcdandhermite}, $Q(\la)=N(\la)G_c(\la)$ for a left invertible $N(\la)$ and a compact square GCRD $G_c(\la)$. Then
$$ P(\la) = \begin{bmatrix}Q(\la)&0
\end{bmatrix}V(\la)^{-1}=N(\la) \begin{bmatrix}G_c(\la)&0 \end{bmatrix} V(\la)^{-1} = N(\la) U(\la)^{-1} G(\la)$$
which proves the statement since $N(\la)U(\la)^{-1}$ is by construction left invertible. \hfill
\end{proof}

The next step is to show how to construct a generic GCRD from a compact GCRD.
\begin{corollary}\label{cor:rankofG}
Let $\{P_i(\la)\}_{i=1}^k$, with $P_i(\la)\in \F[\la]^{m_i\times n}, \;i=1,\ldots,k$,  be a non-empty set of polynomial matrices such that $m:= \sum_{i=1}^k m_i \ge n$, and define their compound matrix $P(\la)$ as in \eqref{eq:compound}. If $G(\la) \in \F[\la]^{\ell \times n}$ is a GCRD of $\{P_i(\la)\}_{i=1}^k$, then $\rank G(\la)=r=\rank P(\la).$
\end{corollary}
\begin{proof}
By Theorem \ref{thm:gcdandhermite}, there exists a compact GCRD $G_c(\la)$ of the given set, and by Proposition \ref{prop:fullrowrank} $\rank G_c(\la)=r$. On the other hand, by Lemma \ref{lem:kernel}, $\ker G(\la)=\ker G_c(\la)$. Thus, by the rank-nullity theorem $\rank G(\la)=n-\dim \ker G(\la)=n-\dim \ker G_c(\la)=r$. \hfill
\end{proof}

\begin{theorem}\label{thm:nonminimalgcd}
Let $P_i(\la)\in \F[\la]^{m_i\times n}, \; i=1,\ldots,k$, be a non-empty set of polynomial matrices such that $m:= \sum_{i=1}^k m_i \ge n$,  let $P(\la)$ be their compound matrix as in \eqref{eq:compound} and assume $r=\rank P(\la) \geq 1$. Then, $G_c(\la) \in \F[\la]^{r \times n}$ is a compact GCRD of $\{ P_i(\la) \}_{i=1}^k$ if and only if $G(\la)$ is a GCRD of $\{ P_i(\la) \}_{i=1}^k$, where 
\begin{equation}\label{eq:gcdgeneral}
G(\la) = Z(\la) \begin{bmatrix}G_c(\la)\\
0\end{bmatrix} \in \F[\la]^{\ell \times n}, \qquad \ell \geq r,
\end{equation}
and $Z(\la) \in \F[\la]^{\ell \times \ell}$ is an arbitrary unimodular matrix, Moreover, any GCRD $G(\la)$ of $\{P_i(\la)\}$ can always be expressed as a (left) linear combination
$$G(\la)=\sum_{i=1}^k L_i(\la)P_i(\la)$$
for some polynomial matrices $L_i(\la) \in \F[\la]^{\ell \times m_i}$; equivalently, there is a polynomial matrix $L(\la)\in\F[\la]^{\ell \times m}$ such that $G(\la)=L(\la)P(\la)$
\end{theorem}
\begin{proof}
Suppose first that $G_c(\la)$ is a compact GCRD of the given set. This implies that $P(\la)=N(\la)G_c(\la)$ for some $N(\la) \in \F[\la]^{m \times r}$. Hence, \[ P(\la) = \begin{bmatrix} N(\la) & 0_{m \times (\ell-r)} \end{bmatrix} Z(\la)^{-1} G(\la) \] and thus $G(\la)$ is a CRD of the same set. On the other hand, $G(\la)=Z_\ell(\la) G_c(\la)$ where $Z_\ell(\la) \in \F[\la]^{\ell \times r}$ contains the leftmost $r$ columns of $Z(\la)$, and hence $G(\la)$ is a GCRD by Lemma \ref{lem:conditionforg}.

For the reverse implication, let us first check that any GCRD $G(\la)$ must have the form \eqref{eq:gcdgeneral}. By definition of GCRD there are matrix polynomials $A(\la) \in \F[\la]^{r \times \ell}, B(\la) \in \F[\la]^{\ell \times r}$ such that $G_c(\la)=A(\la)G(\la)$ and $G(\la)=B(\la)G_c(\la)$. Hence $G_c(\la)=A(\la)B(\la)G_c(\la)$; but $G_c(\la)$ has a right inverse over the field of fractions $\F(\la)$ by Proposition \ref{prop:fullrowrank}, thus $A(\la)B(\la)=I_r$. It follows that $B(\la)$ is left invertible over the ring of polynomials $\F[\la]$ and hence by \cite[Theorem 3.3]{amz} it holds $B(\la)=Z(\la) \begin{bmatrix}I_r\\
0\end{bmatrix}$ for some unimodular $Z(\la)\in\F[\la]^{\ell \times \ell}$, as sought.

Suppose now that $G(\la)$ as in \eqref{eq:gcdgeneral} is a GCRD of the given set. Then, $P(\la)=M(\la)G(\la)$ for some $M(\la) \in \F[\la]^{m \times \ell}$ implying $P(\la)=M(\la) Z_\ell(\la) G_c(\la)$ where $Z_\ell(\la) \in \F[\la]^{\ell \times r}$ contains the leftmost $r$ columns of $Z(\la)$: this shows that $G_c(\la)$ is a compact CRD of $\{ P_i(\la) \}_{i=1}^k$. By Theorem \ref{thm:gcdandhermite}, there exists a compact GCRD of the same set, say, $D_c(\la) \in \F[\la]^{r \times n}$. Hence, by the first part of this proof, $D(\la)=\begin{bmatrix}D_c(\la)\\
0_{(\ell-r)\times n}\end{bmatrix} \in \F[\la]^{\ell \times n}$ is a CRD of $\{ P_i(\la) \}_{i=1}^k$. Thus, since $G(\la)$ is a GCRD, we have  $G(\la)=X(\la)D(\la)$ for some $X(\la) \in \F[\la]^{\ell \times \ell}$. Hence,
\[ G_c(\la) = \begin{bmatrix} I_r & 0 \end{bmatrix} Z(\la)^{-1} X(\la) \begin{bmatrix} I_r\\
0 \end{bmatrix} D_c(\la),\]
proving that $G_c(\la)$ is a GCRD by Lemma \ref{lem:conditionforg}. 

It remains to prove the last part of the statement. Using Theorem \ref{thm:gcdandhermite}, we know that any compact GCRD can be expressed as a left linear combination, say, $G_c(\la)=\sum_{i=1}^k M_i(\la) P_i(\la)$. Therefore, if $G(\la)$ is a GCRD, then for some unimodular $Z(\la)$ and some compact GCRD $G_c(\la)$ it holds
$$G(\la)=Z(\la)  \begin{bmatrix}\sum_{i=1}^k M_i(\la)P_i(\la)\\
0\end{bmatrix} = \sum_{i=1}^k \left( Z(\la) \begin{bmatrix}M_i(\la)\\
0\end{bmatrix}\right) P_i(\la)$$
and the statement follows by defining $L_i(\la):= Z(\la)\begin{bmatrix}M_i(\la)\\
0
\end{bmatrix}$. \hfill
\end{proof}

\begin{remark}
Note that, by applying simultaneously Theorem \ref{thm:kailathisconfusing} and Theorem \ref{thm:nonminimalgcd}, we have that a generic GCRD of $\{ P_i(\la)\}_{i=1}^k$, such that $P(\la)$ as in \eqref{eq:compound} has rank $r\geq 1$, must have the form
\[ G(\la) = U(\la) \begin{bmatrix}D(\la) & 0\\
0 & 0 \end{bmatrix} V(\la)^{-1} \in \F[\la]^{\ell \times n}\] 
where:
\begin{itemize}
\item $U(\la) \in \F[\la]^{\ell \times \ell}$ is an arbitrary unimodular matrix;
\item $V(\la) \in \F[\la]^{n \times n}$ is a unimodular matrix such that $P(\la) V(\la)=\begin{bmatrix}Q(\la) & 0_{m \times (n-r)} \end{bmatrix}$ where $Q(\la) \in \F[\la]^{m \times r}$ has full column rank;
\item $D(\la) \in \F[\la]^{r \times r}$ is a square and compact (and hence regular) GCRD of $\{Q_i(\la)\}_{i=1}^k$, defined as in Theorem \ref{thm:kailathisconfusing}.
\end{itemize}
Satisfying the properties above is, indeed, a sufficient and necessary condition for being a GCRD.
\end{remark}

It is obvious by the definition of a GCRD that, if $G(\la)$ is a GCRD of a given set of polynomial matrices $\{P_i(\la)\}$ and $U(\la)$ is a unimodular matrix, then $U(\la)G(\la)$ is also a GCRD for the same set. In fact, once we have fixed the number of rows of a GCRD, these are the only possible degrees of freedom in determining a GCRD, in the sense that a reverse implication holds. This fact was stated in the literature, but only for the full column rank case; combining it with Theorem \ref{thm:nonminimalgcd}, it characterizes the set of all possible GCRDs.

\begin{theorem}\label{thm:unimodular}
Let $P_i(\la)\in \F[\la]^{m_i\times n}, \; i=1,\ldots,k$, be a non-empty set of polynomial matrices and suppose that $G_1(\la),G_2(\la)$ are two GCRDs of this set such that $G_1(\la)$ and $G_2(\la)$ have the same number of rows. Then $G_1(\la)=U(\la)G_2(\la)$ for some unimodular $U(\la)$.
\end{theorem}
\begin{proof}
Note first that, by definition of GCRD, there are polynomial matrices $A(\la),B(\la)$ such that $G_1(\la)=A(\la)G_2(\la)$ and $G_2(\la)=B(\la)G_1(\la)$.
Suppose first that $G_1(\la)$ and $G_2(\la)$ are compact GCRDs. Then, by Lemma \ref{lem:kernel}, $G_1(\la)$ is right invertible over the field of fractions $\F(\la)$. Hence, $G_1(\la)\!=\!A(\la)B(\la)G_1(\la)$ implies that $I=A(\la)B(\la)$ so that both $A(\la)$ and $B(\la)$ must be unimodular; in particular taking $U(\la):=A(\la)$ proves the statement in this case.

If instead the two GCRDs are not compact, by Theorem \ref{thm:nonminimalgcd} we have that $G_i(\la)$ must have the form
$$ G_i(\la) = Z_i(\la) \begin{bmatrix}R_i(\la)\\
0 \end{bmatrix} $$ for all $i=1,2$ where $Z_1(\la)$ and $Z_2(\la)$ and $V(\la)$ are unimodular whereas $R_1(\la),R_2(\la)$ are compact GCRDs of the same set. By the first part of this proof, there exists a unimodular matrix $V(\la)$ such that $R_1(\la)=V(\la) R_2(\la)$. Thus, 
\[ G_1(\la) = Z_1(\la) \begin{bmatrix}V(\la) & 0\\
0 & I \end{bmatrix} \begin{bmatrix}R_2(\la)\\
0\end{bmatrix} = Z_1(\la) \begin{bmatrix}V(\la) & 0\\
0 & I \end{bmatrix} Z_2(\la)^{-1} G_2(\la), \]
which concludes the proof since $U(\la):=Z_1(\la) (V(\la) \oplus I) Z_2(\la)^{-1}$ is unimodular.  \hfill
\end{proof}

\begin{remark}\label{rem1}
It follows from Theorem \ref{thm:gcdandhermite} that any compact GCRD satisfies 
$P(\la)=N(\la)G_c(\la)$ where $N(\la)$ is left invertible. Hence, any compact GCRD $G_c(\la)$ must contain all the finite eigenvalues of $P(\la)$ since the rank drops of $P(\la)$ must occur in the right factor $G_c(\la)$. Moreover, by Theorem \ref{thm:nonminimalgcd}, every GCRD $G(\la)$ is unimodularly left-equivalent to $\begin{bmatrix}G_c(\la)\\
0\end{bmatrix}$ where $G_c(\la)$ is a compact GCRD. Thus, more generally, any (not necessarily) compact GCRD $G(\la)$ must contain all the finite eigenvalues of $P(\la)$. In addition, Lemma \ref{lem:kernel} shows that $G(\la)$ must also contain all the right minimal indices of $P(\la)$. Any compact GCRD $G_c(\la)$ has the additional property to have full row rank $r$, and therefore it has no left minimal indices. Moreover, the additional degree of freedom of a left unimodular factor $U(\la)$
yields a compact GCRD $U(\la)G_c(\la)$ that can be made row proper,
and therefore has no zeros at infinity in the sense of McMillan \cite{kailath}. (We note that a row proper polynomial matrix may still have infinite eigenvalues in the sense of Gohberg-Lancaster-Rodman \cite{GLR82}.)
\end{remark}

\begin{remark}\label{rem2}
Let $\F \subseteq \K$ be a field extension, and suppose that $D(\la) \in \F[\la]^{p \times n}$ is a CRD (resp. compact CRD) of a set of polynomial matrices $P_i(\la) \in \F[\la]^{m_i \times n}$, $i=1,\dots k$. It is then obvious, by Definition \ref{def:GCRD} and by the invariance of the rank under field extensions, that $D(\la)$ is also a CRD (resp. compact CRD) of the same set of polynomial matrices, when seen as polynomial matrices in $\K[\la]^{m_i \times n}$. Moreover, it follows from Theorem \ref{thm:gcdandhermite} and Theorem \ref{thm:nonminimalgcd} that if $G(\la) \in \F[\la]^{\ell \times n}$ is a GCRD (resp. compact GCRD) over $\F[\la]$, it is also a GCRD (resp. compact GCRD) over $\K[\la]$.

For example, let us consider the case $\F=\R$ and $\K=\C$. The theory we developed shows that every non-empty set of real polynomial matrices $P_i(\la) \in \R[\la]^{m_i \times n}$, $i=1,\dots,k$ has a real GCRD $G(\la) \in \R[\la]^{\ell \times n}$ for all $\ell \geq r$ where $r=\rank P(\la)$ in \eqref{eq:compound}; the observations above additionally show that such $G(\la)$ is also a complex GCRD, i.e., a GCRD for the same $P_i(\la)$ seen as complex polynomial matrices. Moreover, by Theorem \ref{thm:unimodular}, any other real GCRD with the same number of rows has the form $U(\la)G(\la)$ where $U(\la) \in \R[\la]^{\ell \times \ell}$ is a real unimodular matrix; and any other \emph{complex} GCRD with the same number of rows has the form $V(\la)G(\la)$ where $V(\la) \in \C[\la]^{\ell \times \ell}$ is a complex unimodular matrix. As a consequence, if one has computed a complex nonreal GCRD for a given set of real polynomial matrices, it is always possible to obtain a real GCRD for the same set by some appropriate left unimodular multiplication.
\end{remark}

\section{State-space realizations and feedback} \label{sec:feedback}
In this section and the following ones, we restrict our attention to the case where the base field $\F$ is the field of complex numbers $\C$, with the goal of developing and testing a novel numerical algorithm. The case $\F=\R$ is analogous, but one has to rely on the real staircase form; for simplicity of exposition, we just focus on the complex case.

Before describing the new algorithm for GCRD extraction of a polynomial matrix, we briefly recall the properties of (generalized and classical) state-space realizations and the effect of feedback on them. Since only proper transfer functions have a standard state-space realization, 
one either has to make use of a generalized state-space or pencil realization, or perform a change of variables $\la = 1/\mu$ to 
then realize $P(1/\mu)$ in standard state-space form. In \cite{emre,SilV81} an approach to computing GCRDs based on state-space realization and feedback was described, but there an assumption was made that \eqref{eq:compound} has full rank. Those methods were based on the notion of $(A,B)$-invariant subspaces \cite{Wonham} which required the use of state-space models. In this section, we first recall those ideas, using the method described in \cite{SilV81}.

Our starting point is an expansion of the compound matrix \eqref{eq:compound} in the monomial basis $P(\la)=: P_0 + P_1 \la + \ldots + P_d \la^d,$ where $d$ is the degree of \eqref{eq:compound}. A pencil realization of $P(\la)$ is a quintuple of matrices ${A,B,C,D,E}$ such that the system matrix 
of $P(\la)$ is given by
\[
S_\la(\lambda) := \left[ \begin{array}{c|c} A-\lambda E & \la B \\
\hline C & D \end{array} \right] , \quad P(\la)= C(\la E-A)^{-1}\la B+D.
\]
Note that $P(\lambda)$
is the Schur complement of its system matrix $S_\la(\la)$, with respect to the (invertible) leading submatrix $A-\la E$. A classical solution for this is given by the companion-like form 
\begin{equation} \label{poly}
S_\la(\la)= 
\left[ \begin{array}{ccccc|c} I_n & -\la I_n & & & \\  &  I_n & -\la I_n & & \\  & & \ddots & \ddots & \\
& & & I_n & -\la I_n & \\
& & & & I_n & -\la I_n
\\ \hline P_d & \ldots & \ldots & P_2 & P_1 & P_0
\end{array} \right].
\end{equation}
The pencil realization \eqref{poly} is also {\it irreducible} since neither the submatrix
$ 
\left[ \begin{array}{c} A-\lambda E  \\
\hline C \end{array} \right]$
 nor the submatrix 
$\left[ \begin{array}{c|c} A-\lambda E & \la B
\end{array} \right]$
have any finite Smith zeros \cite{Rosenbrock,vvk}.

\begin{remark}
It is well known that there exist unimodular transformations $U(\la)$ and $V(\la)$ such that 
$$  U(\la)S_\la(\la)V(\la) = \left[ \begin{array}{c|c} I_{dn} & 
\\ \hline  & P(\la) \end{array} \right].
$$
This shows clearly that the Smith zeros of $P(\la)$
are also those of $S_\la(\la)$ and that the pencil $S_\la(\la)$ has normal rank $dn+r$. In other words, $S_\la(\la)$ is a linearization of $P(\la)$ \cite{GLR82} and it can be used to compute the finite zeros of $P(\la)$ via the staircase algorithm, followed by the QZ algorithm.
\end{remark}

\bigskip

If we use the change of variables $\la=1/\mu$, then a standard state space realization of 
$$  P(1/\mu)= P_0 + P_1/\mu + \ldots + P_d/\mu^d
$$
is given by the following system matrix
\begin{equation} \label{polymu}
S_\mu(\mu):= \left[ \begin{array}{c|c} \mu I_{dn}-A & -B 
\\ \hline C & D
\end{array} \right] :=
\left[ \begin{array}{ccccc|c} \mu I_n & - I_n & & & \\  &  \mu I_n & - I_n & & \\  & & \ddots & \ddots & \\
& & & \mu I_n & -I_n & \\
& & & & \mu I_n & - I_n
\\ \hline P_d & \ldots & \ldots & P_2 & P_1 & P_0
\end{array} \right],
\end{equation}
which is not observable if $\rank P_d$ is smaller than $n$, but it is controllable. We will see that controllability suffices for our purpose. The relation 
 $\label{mula} S_\la(\la) = (\la I_{dn} \oplus I_m) S_\mu(\mu)$
also shows that \eqref{poly} and \eqref{polymu} have the same normal rank.
If we assume that $P_0$ has full column rank, then $P(\la)$ has full normal rank $n$, $P(\la)$ and $S_\la(\la)$ have no zeros at $\la=0$ and $P(1/\mu)$ and $S_\mu(\mu)$ have no zeros at $\mu=\infty$. 

\bigskip

Theorem \ref{th:decomp} below describes a decomposition that holds for the system matrix $S_\mu(\mu)$, under the assumption that $P_0$ has full rank $n$.
This approach is essentially described in the papers \cite{emre,SilV81}, but the present derivation avoids elaborate system theoretic concepts. The full rank assumption implies that the method proposed in \cite{emre,SilV81} can only be applied to the case where \eqref{eq:compound} has full column rank (up to a linear change of variable $\la \mapsto \la-\alpha$, this can be assumed to be equivalent to a full rank $P_0$). The main advantage of the standard state space model for $P(1/\mu)$ is that 
it has only four matrices $\{A,B,C,D\}$ and that subsequent transformations of the model are simpler. It also shows that the construction of a particular feedback solves the GCRD extraction problem.

\begin{theorem} \label{th:decomp}
Let $P_0$ have full column rank $n$ and let $ \left[ \begin{array}{c|c} \mu I_{dn}-A & -B 
\\ \hline C & D
\end{array} \right]$ be as in \eqref{polymu}. Then
there exist a unitary transformation $U\in \C^{dn\times dn}$ and a ``feedback" matrix  $F\in \C^{n\times dn}$
such that 
\begin{equation} \label{decomp}
 \left[ \begin{array}{cc|c}
\mu I_{d_f}-\hat A_{11} & - \tilde A_{12}  & - B_1  \\
0 & \mu I_{d_0} - \tilde A_{22} & - B_2 \\ \hline  0 & \tilde C_2 & D \end{array} \right]:=
  \left[ \begin{array}{c|c}
 U^* & 0  \\ \hline
0 & I_m  \end{array} \right] S_\mu(\mu) \left[ \begin{array}{c|c}
 U & 0  \\ \hline
F & I_n  \end{array} \right] 
\end{equation} 
where $\hat A_{11}$ is square and contains the $d_f$ finite eigenvalues of the pencil $S_\mu(\mu)$, $\tilde A_{22}$ is square and has its $d_0=dn-d_f$ eigenvalues at 0 and the pair $(\tilde A_{22},(I_m-P_0P_0^\dagger)\tilde C_2)$ is observable.
\end{theorem}
\begin{proof}
The proof is constructive. Apply the feedback 
$$   \hat F := \left[\begin{array}{cccc} \hat F_d & \ldots & \hat F_1  \end{array}\right] := - P_0^{\dagger}\left[\begin{array}{cccc} P_d & \ldots &  P_1 \end{array}\right] = - P_0^{\dagger}C
$$
and define the matrices $\hat A$ and $\hat C$ as follows
$$  \hat A := A+B\hat F, \quad \hat C := C+D\hat F=  (I_m-P_0P_0^{\dagger}) \left[\begin{array}{cccc} P_d & \ldots & P_1  \end{array}\right] =  (I_m-P_0P_0^{\dagger}) C $$
to obtain the new system matrix 
$$  \left[ \begin{array}{c|c}
 \mu I_{dn} - \hat A  & - B \\ \hline  \hat C &  D \end{array} \right] :=
  S_\mu(\mu) \left[ \begin{array}{c|c}
 I_{dn} & 0  \\ \hline
\hat F & I_n  \end{array} \right]. $$
Now compute the orthogonal state-space transformation $U$ that displays the $d_f$ unobservable modes of the pair $(\hat A, \hat C)$ in a submatrix $\hat A_{11}$, as follows \cite{Van81}:
$$
 \left[ \begin{array}{cc|c}
 \mu I_{d_f} - \hat A_{11}  & - \hat A_{12}  & - B_1 \\
 0  & \mu I_{d_0} - \hat A_{22}  & -  B_2 \\
 \hline 0 & \hat C_2 & D  \end{array} \right] :=  \left[ \begin{array}{c|c}
U^* & 0  \\ \hline
0 & I_{m}  \end{array} \right]
\left[ \begin{array}{c|c}
 \mu I_{dn} - \hat A  & - B \\ \hline \hat C & D \end{array} \right]
 \left[ \begin{array}{c|c}
 U & 0  \\ \hline
0 & I_n  \end{array} \right] .
$$
Finally, apply an additional feedback law $\tilde F_2$ to place all the eigenvalues of the subsystem 
$(\hat A_{22}, B_2)$ at zero (which is known as deadbeat control \cite{Van84}). This is possible since that subsystem inherits its controllability from that of the full system.
The resulting matrix $\tilde A_{22}:=\hat A_{22}+\hat B_2\tilde F_2$ is then nilpotent
and its index of nilpotency can be chosen to be minimal~\cite{Van84}:
$$
 \left[ \begin{array}{cc|c}
 \mu I_{d_f} - \hat A_{11}  & - \tilde A_{12}  & - B_1 \\
 0  & \mu I_{d_0} - \tilde A_{22}  & -  B_2 \\ \hline 
 0 & \tilde C_2 & D  \end{array} \right] :=   \left[ \begin{array}{cc|c}
 \mu I_{d_f} - \hat A_{11}  & - \hat A_{12}  & - B_1 \\
 0  & \mu I_{d_0} - \hat A_{22}  & -  B_2 \\
 \hline 0 & \hat C_2 & D  \end{array} \right].
 \left[ \begin{array}{cc|c}
 I_{d_f} & 0  & 0 \\ 0 & I_{d_0} & 0 \\\hline
0 & \tilde F_2 & I_n  \end{array} \right].
$$
Combining the different steps and setting
 $F:=\hat FU + \left[\; 0 , \tilde F_2 \right]$ yields the desired result \eqref{decomp}. \hfill
\end{proof}

Notice that the deadbeat control part can be avoided, because before the feedback $\hat F$ was applied, all eigenvalues of $A$ were at $\mu=0$. It suffices to only apply the part of the feedback that cancels the finite zeros and leaves the other eigenvalues unchanged.
How this can be done is shown in \cite{SilV81}.

\begin{corollary}
It follows from Theorem \ref{th:decomp} that  $P(\la)= N(\la)G(\la)$  with
$$ N(\la)=  D  + \tilde C_2 (I_{d_0}-\la \tilde A_{22})^{-1} \la B_2, \quad \mathrm{and} \quad G(\la)=I_n -F(I_{dn}-\la A)^{-1} \la B
$$ 
and that $N(\la)$ is left invertible and $G(\la)$ is a GCRD for the blocks of the compound matrix $P(\la)$.
\end{corollary}
\begin{proof}
It follows from \eqref{mula} and \eqref{decomp} that we also have
$$  \left[ \begin{array}{c|c}
I_{dn}- \la U^*AU  & -\la U^*B  \\ \hline
CU & D  \end{array} \right] \left[ \begin{array}{c|c}
 I_{dn} & 0  \\ \hline
F & I_n  \end{array} \right] = \left[ \begin{array}{cc|c}
I_{d_f}- \la \hat A_{11}  & -\la \tilde A_{12}  & -\la B_1  \\
0 & I_{d_0}-\la \tilde A_{22}  & -\la B_2 \\ \hline  0 & \tilde C_2 &  D \end{array} \right].
$$
The transfer function of the system quadruple on the right hand side is 
$$  D  + \tilde C_2  (I_{d_0}-\la \tilde A_{22})^{-1}\la B_2$$
which is polynomial since $\tilde A_{22}$ is nilpotent, and has no finite zeros since there exists a unitary matrix 
$V$ such that 
$$  \left[ \begin{array}{c|c}
I_{d_0} & 0 \\ \hline 0 & V \end{array} \right] 
 \left[ \begin{array}{c|c}
I_{d_0}- \la \tilde A_{22} & -\la B_2 \\ \hline (I_m-P_0P_0^\dagger)C &  P_0 \end{array} \right] =   \left[ \begin{array}{c|c}
I_{d_0}- \la \tilde A_{22} & -\la B_2 \\ \hline  \tilde C_V & 0  \\  0 & \tilde P_0 \end{array} \right] 
$$
where $\tilde P_0$ is invertible, and $(\tilde A_{22},\tilde C_V)$ is observable. Therefore that pencil has 
no finite zeros. As shown in Lemma \ref{lem:feedback}, the feedback operation can be written as the matrix factorization $P(\la)=N(\la)G(\la)$
where the polynomial matrices $P(\la)$ and $G(\la)$ are realized by the system matrices 
$$   \left[ \begin{array}{c|c}
I_{0} -\la \tilde A_{22} & -\la B_2  \\ \hline
\tilde C_2 & D  \end{array} \right], \quad \mathrm{and} \quad
 \left[ \begin{array}{c|c}
I_{dn}-\la A & -\la B  \\ \hline
-F & I_n  \end{array} \right]
$$
respectively.
Moreover, the coefficients $G_i, i=0,\ldots,d$ and $N_i, i=1,\ldots,d$ of the polynomial expansions $G(\la)=\sum_{i=0}^d G_i\la^i$ and $N(\la)=\sum_{i=0}^d N_i\la^i$ are given by
$$ G_0=I_n,\; G_i=-F_i, \; i=1,\ldots,d, \quad N_0=P_0, \; N_i=\tilde C_2 \tilde A_{22}^{i-1}B_2,  \; i=1,\ldots,d.
$$
 \hfill
\end{proof}

\section{Generalized state-space approach}  \label{sec:gssm}
In this section we describe a novel variation of the approach of Section \ref{sec:feedback}, leading to a new algorithm. Let us consider a general $m\times n$ polynomial matrix $P(\la)$ of normal rank $r$ and represent it as the transfer function of a generalized state space model 
\begin{equation} \label{gssm}
S_P(\la)= 
\left[ \begin{array}{c|c} A-\la E & B  \\ \hline C & D
\end{array} \right] :=
\left[ \begin{array}{ccccc|c} I_n & -\la I_n & & & \\  &  I_n & -\la I_n & & \\  & & \ddots & \ddots & \\
& & & I_n & -\la I_n & \\
& & & & I_n & - I_n
\\ \hline P_d & \ldots & \ldots & P_1 & P_0 & 0
\end{array} \right]
\end{equation}
which has state-space dimension $(d+1)n$ and is irreducible for all finite points. This pencil has normal 
column rank equal to $(d+1)n+r$, and its finite zeros are the finite eigenvalues of the pencil $S_P(\la)$.
Moreover, if $\Lambda(\la)$ is the Smith form of the matrix $P(\la)$, then the Smith form of $S_P(\la)$ defined in \eqref{gssm}, and $S_\la(\la)$
given in \eqref{poly} are, respectively,  
$$  I_{(d+1)n} \oplus \Lambda(\la) \quad\mathrm{and} \quad I_{dn} \oplus \Lambda(\la).
$$
It is also easy to see that the pencil $S_P(\la)$
is strictly equivalent to $I_n \oplus S_\la(\la)$
which implies that the left and right minimal indices and the finite and (non-trivial) infinite zeros of \eqref{poly} and \eqref{gssm} are the same.

The following decomposition can then be obtained by running the staircase reduction on the $(dn+m)\times (dn+n)$ subpencil $S_\la(\la)$, given in \eqref{poly}, of the system matrix $S_P(\la)$.

\begin{theorem} \label{th:transform}
Let the $m\times n$ polynomial matrix $P(\la)$ have normal rank $r$, then there exist 
unitary transformations $Q\in \C^{dn\times dn}$ and $Z\in \C^{(d+1)n\times (d+1)n}$ satisfying
\[
\left[ \begin{array}{c|c} Q^* & \\ \hline &  I_m 
\end{array} \right] \!
\left[ \begin{array}{ccccc} I_n & -\la I_n & & \\  &  I_n & -\la I_n &  \\  & & \ddots & \ddots  \\
& & & I_n & -\la I_n  
\\ \hline P_d & \ldots & \ldots & P_1 & P_0 
\end{array} \right] \! Z
\]
$$=\left[ \begin{array}{cccc} A_{11}-\la E_{11} & A_{12}-\la E_{12} & A_{13}-\la E_{13}  & A_{14}-\la E_{14} \\
0 & A_{22}-\la E_{22}  & A_{23}-\la E_{23} & A_{24}-\la E_{24} \\ 
0 & 0 & A_{33}-\la E_{33} & A_{34}-\la E_{34} \\ 
\hline  0 & 0 & 0 & C_4
\end{array} \right],
$$
where:
\begin{itemize}
\item the pencil $A_{22}-\la E_{22}$ is regular and contains the finite zeros of  $S_\la(\la)$;
    \item the pencil $A_{11}-\la E_{11}$ contains
the $n-r$ right Kronecker blocks of
$S_\la(\la)$;
\item the pencil
$$\left[ \begin{array}{cc}  
 A_{33}-\la E_{33}  & A_{34}-\la E_{34} \\ 
\hline   0 & C_4
\end{array} \right],$$
contains the infinite zeros and the $m-r$ left Kronecker blocks of $S_\la(\la)$;
\item the submatrix 
$C_4$ has full column rank.
\end{itemize}
\end{theorem}
\begin{proof} The staircase form, described in Lemma \ref{lem:staircase}, precisely separates the right minimal indices and the finite zeros of a given pencil, from its remaining structural elements, which in our case are the infinite zeros and the left minimal indices of the pencil. It does so using a unitary equivalence transformation that separates the different groups of structural elements in a block upper triangular form. In a first stage, the right Kronecker blocks and finite zeros are separated from the infinite zeros and left Kronecker blocks. In second stage an updating unitary equivalence transformation separates further the right Kronecker blocks from the finite zeros, to generate a three-way decomposition, as described in Lemma \ref{lem:staircase}. The first step of stage 1, is a row transformation compressing the rows of the coefficient of $\la$. But the coefficient of $\la$ in the pencil $S_\la(\la)$ is already in row-compressed form since its top part is the submatrix $I_{dn}$. The next step then compresses the columns of the submatrix 
$$  \left[\begin{array}{cccc}P_d  & \ldots & P_1 & P_0  \end{array} \right] 
$$
to the right, 
to produce the matrix $C_4$ of full column rank $n_4$. The algorithm then proceeds further with the leading  subpencil of dimension $dn\times(dn-n_4)$, and applies a similar sequence of row and columns compressions to subpencils of decreasing dimensions.
The special block structure of the left transformation $Q^*\oplus I_m$  thus follows from the fact that the very first step of the staircase algorithm can be skipped here.  \hfill
\end{proof}

In Theorem \ref{th:rankembed} below we analyze how the decomposition of Theorem \ref{th:transform} can be combined with feedback to yield 
a transfer function with the desired properties.

\begin{theorem} \label{th:rankembed}
Let the $m\times n$ polynomial matrix $P(\la)$ have normal rank $r$, then there exists a feedback matrix $F$ and unitary matrices $Q$ and $Z$ such that 
\[
\left[ \begin{array}{cc|c} Q^* & & \\  & I_n &  \\ \hline & &  I_m 
\end{array} \right] 
\left[ \begin{array}{ccccc|c} I_n & -\la I_n & & & \\  &  I_n & -\la I_n & & \\  & & \ddots & \ddots & \\
& & & I_n & -\la I_n & \\ & & & & I_n & - I_n   
\\ \hline P_d & \ldots & \ldots & P_1 & P_0 & 0 
\end{array} \right]  \left[ \begin{array}{c|c} Z &  \\ \hline FZ & I_n  
\end{array} \right] 
\]
\begin{equation}\label{eq:FQZ} 
=\left[ \begin{array}{cccc|cc} A_{11}-\la E_{11} & A_{12}-\la E_{12} & A_{13}-\la E_{13}  & A_{14}-\la E_{14} & 0 & 0 \\ 
0 & A_{22}-\la E_{22}  & A_{23}-\la E_{23} & A_{24}-\la E_{24} & 0 & 0 \\ 
0 & 0 & A_{33}-\la E_{33} & A_{34}-\la E_{34} & 0 & 0 \\
0 & 0 & Z_3 & Z_4 & -I_{r} & 0 \\
Z_1 & 0 & 0 & 0 & 0 & -I_{n-r} \\
\hline  0 & 0 & 0 & C_4 & 0 & 0 
\end{array} \right],
\end{equation}
where $A_{22}-\la E_{22}$ is square, contains the finite zeros of $P(\la)$ and is unobservable in this realization, 
the pencil $\left[\begin{array}{cc} A_{11}-\la E_{11}\\Z_1  \end{array}\right]$ is unimodular and 
unobservable in this realization,
and the subpencil
\begin{equation} \label{embed} 
\left[ \begin{array}{cc} 
A_{33}-\la E_{33} & A_{34}-\la E_{34} \\ 
Z_3 & Z_4 \end{array} \right],
\end{equation}
is unimodular. 
\end{theorem}
\begin{proof}
In order to obtain the decomposition \eqref{eq:FQZ} we first apply a feedback $F=\left[ \begin{array}{cccc} 0 & \ldots & 0 & I_n \end{array} \right]$, to eliminate the coupling term $I_n$ in the last block row of $(A-\la E)$~:
\[\left[ \begin{array}{cccc|c} I_n & -\la I_n & & & \\  & \ddots & \ddots & \\ & & I_n & -\la I_n & \\ & & & 0 & - I_n   
\\ \hline P_d & \ldots & P_1 & P_0 & 0 
\end{array} \right] = \left[ \begin{array}{cccc|c} I_n & -\la I_n & & & \\  & \ddots & \ddots & \\ & & I_n & -\la I_n & \\ & & & I_n & - I_n   
\\ \hline P_d & \ldots & P_1 & P_0 & 0 
\end{array} \right] \left[ \begin{array}{c|c} I_{(d+1)n} &  \\ \hline F & I_n  
\end{array} \right].
\]
Then we integrate in this decomposition the transformations $Q$ and $Z$
of Theorem \ref{th:transform}. It follows from the construction of that Theorem, that the pencils
$A_{11}-\la E_{11}$ and $\left[ \begin{array}{cc} 
A_{33}-\la E_{33} & A_{34}-\la E_{34} \end{array} \right]$ have full row rank for all finite $\la$, and that their leading coefficient has full row rank as well. By \cite[Lemma 5.2]{nvdpen}, these pencils are then right invertible and, as a consequence \cite[Theorem 3.3]{amz}, they always have a unimodular embedding \eqref{embed}, realized by the matrices $Z_1$ and  $\left[ \begin{array}{ccc}
Z_{3}  &  Z_{4} \end{array} \right]$, respectively. This matrices $Z_1$, $Z_3$ and $Z_4$ are then placed in the bottom block row of $(A-\la E)$ by the updated  feedback matrix 
$$  F =  \left[ \begin{array}{cccc} 0 & 0 & 0 & I_n \end{array} \right] -\left[ \begin{array}{cccc} 0 & 0 & Z_3 & Z_4 \\ Z_1 & 0 & 0 & 0  \\ \end{array} \right]Z^T.
$$
A construction for the embedding matrices $Z_1$ and $\left[ \begin{array}{cccc} Z_3 & Z_4 \end{array} \right]$ is e.g. given in \cite{BeelenV}. 
\hfill
\end{proof}

We have now all the ingredients to construct a factorization, as stated in the following theorem.

\begin{theorem} \label{th:rankfactor}
The decomposition of Theorem \ref{th:transform} of the polynomial matrix $P(\la)$ of normal rank $r$, defined by its system matrix $S_P(\la)$ described in \eqref{gssm}, yields the polynomial factorization $P(\la)=N_r(\la)G_c(\la)$ where $N_r(\la)\in \C[\la]^{m\times r}$ has full rank $r$ for all finite $\la$ and $G_c(\la)\in \C[\la]^{r\times n}$ contains all the finite zeros of $P(\la)$. Moreover, $N_r(\la)$ and $G_c(\la)$ are given by their system matrices 
\begin{equation} \label{eq:lowrank}
S_{N_r}(\la):=\left[ \begin{array}{cc|c}
A_{33}-\la E_{33}  & A_{34}-\la E_{34} & 0 \\ 
Z_3 & Z_4 & -I_r \\
\hline  0 & C_4 & 0 
\end{array} \right], \; S_{G_c}(\la):=
\left[ \begin{array}{c|c} A -\la E & B  \\ \hline -[I_r, 0]F & [I_r, 0]
\end{array} \right].
\end{equation}
\end{theorem}
\begin{proof}
We obtain from the above feedback construction and from Lemma \ref{lem:feedback} the factorization $P(\la)=N(\la)G(\la)$ where the two factors are given by the system matrices 
$$ S_{N}(\la):=\left[ \begin{array}{cc|cc}
 A_{33}-\la E_{33} & A_{34}-\la E_{34} & 0 & 0 \\ 
  Z_3 & Z_4 & -I_r & 0\\
\hline  0 & C_4 & 0 & 0
\end{array} \right], \; S_{G}(\la):=
\left[ \begin{array}{c|c} A -\la E & B  \\ \hline -F & I_n
\end{array} \right].
$$
The matrix $N(\la)$ is clearly zero in its last $n-r$ columns and can thus be written as $N(\la)=N_r(\la)\left[\begin{array}{cc} I_r & 0\end{array}\right]$. And this right factor 
$\left[\begin{array}{cc} I_r & 0\end{array}\right]$ can be absorbed in the product $G_c(\la)=\left[\begin{array}{cc} I_r & 0\end{array}\right]G(\la)$, where the system matrices of 
$N_r(\la)$ and $G_c(\la)$ are as in \eqref{eq:lowrank}.
The two factors are polynomial because their system matrices have unimodular matrices as state-transition pencil. The constructed factor $N_r(\la)$ has full column rank $r$ since its
system matrix $S_{N_r}(\la)$ is strictly equivalent to the pencil
$$   I_r\oplus \left[ \begin{array}{cc}
A_{33}-\la E_{33}  & A_{34}-\la E_{34} \\ 
\hline  0 & C_4 
\end{array} \right],
$$
which has only left minimal indices and infinite zeros.
\hfill
\end{proof}

The computation of the coefficients $N_i$ and $G_i$ of the expansions $N_r(\la)=\sum_{i=0}^d N_i\la^i$ and  $G_c(\la)=\sum_{i=0}^d G_i\la^i$ are now given in terms of the matrices 
$$\left[\begin{array}{cccc} \hat F_d & \ldots & \hat F_1 & \hat F_0\end{array}\right] := \hat F:=\left[\begin{array}{cc} I_r & 0\end{array}\right]F
$$
and
$$    \left[ \begin{array}{c|c}
\hat A - \la \hat E  & -\hat B \\ 
\hline  \hat C & 0
\end{array} \right] :=  \left[ \begin{array}{cc|c}
A_{33}-\la E_{33}  & A_{34}-\la E_{34} & 0 \\ Z_3 & Z_4 & -I_r \\ 
\hline  0 & C_4 & 0
\end{array} \right]
$$
as follows
\[
G_0=I_n-\hat F_0,\; G_i=-\hat F_i, \; i=1,\ldots,d, \quad N_i=\hat C (\hat A^{-1}\hat E)^{i}\hat A^{-1}\hat B,  \; i=0,\ldots,d.
\]

\begin{remark} \label{rem:boundG}
The generalized state space approach has the disadvantage to require the computation of two large transformation matrices, but its main advantage is that the norm of the feedback matrix $F$ can be
chosen to be bounded. All calculations of the decomposition are completely based on orthogonal transformations in the real case, and unitary transformations in the complex case. The final construction of the GCRD is directly obtained from the feedback matrix $F$. Typically, its submatrices $Z_i$ are chosen to have orthonormal rows (see \cite{BeelenV}) and it then follows that each row of the matrix $G_c(\la)$ has Frobenius norm 1. We point out that the construction of the factor $G_c(\la)$ only requires the product $\left[\begin{array}{cc} I_r & 0\end{array}\right]F$. This implies that the matrix $Z_1$ does not need to be computed and that the separation of the 
blocks $(A_{11}-\la E_{11})$ and $(A_{22}-\la E_{22})$ in stage 2 of the 
staircase algorithm does not need to be performed.
\end{remark}

\begin{remark} Note that we can construct from $N_r(\la)$ and $G_c(\la)$ also factorizations $P(\la)=N_\ell(\la) G_\ell(\la)$ 
where $N_\ell(\la) \in \C[\la]^{m\times \ell}$ and  $G_\ell(\la) \in \C[\la]^{\ell \times n}$, for any $r < \ell \leq n$ by appending columns to $N_0$ that are linear independent from those of $N_0$ and appending zero rows to $G_c(\la)$.
\end{remark}

\subsection{Appendix to Section \ref{sec:gssm}: Technical results}
We recall here some basic lemmata about system pencils.

\begin{lemma}  \label{lem:feedback} The following block transformation of the system matrix of $P(\la)$
$$  S_{N}(\la)= \left[ \begin{array}{c|c} A+BF-\la E & B  \\ \hline C+DF & D
\end{array} \right]:= S_P(\la)\left[ \begin{array}{c|c} I & 0  \\ \hline F & I_n
\end{array} \right].
$$
extracts from the transfer function $P(\la)= N(\la) G(\la)$ a right factor $G(\la)$ defined by the system matrix 
$$  S_G(\la)= \left[ \begin{array}{c|c} A-\la E & B  \\ \hline -F & I_n
\end{array} \right]
$$
provided the pencil $(A+BF-\la E)$ is regular.
\end{lemma}
\begin{proof}
It is well known that the system matrix of a product can be obtained from a product of (expanded) system matrices
(see \cite{Van90}). Applying this to the system matrices $S_N(\la)$ and $S_G(\la)$, we obtain the system matrix for the product $N(\la)G(\la)$ as follows~: 
$$   \left[ \begin{array}{cc|c} A+BF-\la E & -BF &  B  \\ 0 & A-\la E & B \\ \hline C+DF & -DF & D
\end{array} \right] =  \left[ \begin{array}{cc|c} A+BF-\la E & 0 &  B  \\ 0  & I & 0 \\ \hline C+DF & 0 & D
\end{array} \right]  \left[ \begin{array}{cc|c} I & 0 & 0 \\ 0 & A-\la E & B  \\ \hline  0 & -F & I_n
\end{array} \right].
$$
This is a valid Rosenbrock system quadruple, provided the pencil $(A+BF-\la E)$ is regular.
Using a trivial block row and block column transformation, we then find that the left system matrix is strict system equivalent to 
$$   \left[ \begin{array}{cc|c} A+BF-\la E & 0 &  0  \\ 0 & A-\la E & B \\ \hline C+DF & C & D
\end{array} \right] \sim  \left[ \begin{array}{c|c} A-\la E &   B  \\ \hline C & D
\end{array} \right] .
$$
\hfill
\end{proof}

We also recall basic properties of the staircase form of pencils used in the paper. For details of the proof we refer to \cite{vd79}, where a lower triangular version is given.
\begin{lemma}[Staircase form] \label{lem:staircase}
Let $A-\la E\in \C[\la]^{m\times n}$ be a general $m\times n$ matrix pencil of normal rank $r$. Then there exist unitary  transformations $Q\in \C^{m\times m}$ and $Z\in \C^{n\times n}$ such that 
$$
Q^*(A-\la E)Z =\left[ \begin{array}{ccc} A_{11}-\la E_{11} & A_{12}-\la E_{12} & A_{13}-\la E_{13} \\ 
0 & A_{22}-\la E_{22}  & A_{23}-\la E_{23}  \\ 
0 & 0 & A_{33}-\la E_{33} \\ 
\end{array} \right],
$$
where $A_{11}-\la E_{11}$ contains
the $n-r$ right Kronecker blocks of  $A-\la E$,  
$A_{22}-\la E_{22}$ is regular and contains the finite zeros of  $A-\la E$,
and $A_{33}-\la E_{33}$ contains the infinite zeros and $m-r$ left Kronecker blocks of $A-\la E$.
\end{lemma}

\section{Numerical aspects}  \label{sec:numerics}
Numerical algorithms for the computation of a GCRD of a polynomial matrix $P(\la)$ can be classified in four categories. The classical (algebraic) approach uses unimodular transformations to construct the Hermite form of the matrix $P(\la)$ (see \cite{Friedland}), from which the GCRD can be derived. The worst case complexity of this approach is significantly larger than for other methods and it is easy to see that coefficient growth can be quite dramatic since numerical pivoting techniques are precluded. A positive aspect of this approach is that it also applies to matrices with deficient normal rank.
A second class of methods is based on spectral divisors of the submatrices $P_i(\la)$ forming $P(\la)$ \cite{GLR82}, but these methods apply only to square submatrices $P_i(\la)$ and also require the computation of their zeros.
A third class of methods is related to resultant matrices \cite{Bitmead,moness,Henrion}.
These matrices and their left null vectors can be used to characterize the left inverses of $P(\la)$ and from there also construct the GCRD. The methods using this approach perform a numerical elimination of linear dependent rows of the Sylvester matrix, but the elimination procedure and corresponding echelon form can be seen to make implicit use of unimodular transformations. This again precludes standard pivoting techniques that are important 
for numerical stability. Moreover, it is unclear how (or if) these methods can be adapted to the case where \eqref{eq:compound} is not full rank. A fourth class of methods is based on state-space realizations \cite{emre,SilV81},
as explained in Section \ref{sec:feedback}.

The method proposed in this paper does not belong to any of these classes, but it is close to the fourth class. Our algorithm\footnote{A MATLAB implementation is freely avaliable from github at the link \href{https://github.com/VanDoorenPaul/Greatest-Common-Divisor}{https://github.com/VanDoorenPaul/Greatest-Common-Divisor}.} uses a \emph{generalized} state-space realization of the polynomial matrix, its staircase reduction to find its decoupling part and state feedback to implement the factorization defining a GCRD. 
Its complexity is that of the staircase algorithm, which is cubic in the state-space dimension (which is smaller than the size of the Sylvester matrix used in the third class of methods). This approach has the additional advantage that all rank decisions are based on orthogonal state space transformations only, but it does not take away the inherent ill-condition of the problem.

We now describe some numerical experiments that we performed to test the new algorithm. We first considered the example 
$$ P(\la)=\left[\begin{array}{rr} 1 & 1 \\  1 & 0 \\ 5 & 2 \\-1 & -1 \end{array}\right] +
\left[\begin{array}{rr} 2 & 0 \\  2 & 2 \\ 3 & 4 \\ 1 & 1 \end{array}\right] \la +
\left[\begin{array}{rr} 0 & 1 \\  1 & 1 \\ 2 & 0 \\ 1 & 1 \end{array}\right] \la^2 +
\left[\begin{array}{rr} 0 & 0 \\  0 & 0 \\ 0 & 1 \\ 0 & 0 \end{array}\right] \la^3,
$$
drawn from \cite{Bitmead}, and for which the authors of \cite{Bitmead} computed a compact GCRD 
$$ G(\la)=\left[\begin{array}{rr} 5 & 2 \\ 1 & 0 \end{array}\right] +
\left[\begin{array}{rr} 2 & 3 \\  0 & 1  \end{array}\right] \la .
$$
Our algorithm uses only orthogonal transformations on the system matrix pencil of $P(\la)$. Therefore, if we normalize $P(\la)$ to
have a norm equal to $1$, all rank decisions will be on submatrices bounded in norm by $1$. For this example, we obtained the computed factorization
$$ N_r(\la)\approx
\left[\begin{array}{rr} 1.4633 & -0.9267 \\ 1.1785 &  0.7818 \\ 6.4620 &  0.4920 \\ -1.4633 &  0.9267 \end{array}\right] +
\left[\begin{array}{rr}  1.1785 &  0.7818 \\ 2.6418 & -0.1449 \\ 1.4633 & -0.9267 \\ 2.6418 & -0.1449 \end{array}\right] \la +
\left[\begin{array}{rr} 0.0000 &  -0.0000 \\ 0.0000 &  -0.0000 \\ 1.1785 &   0.7818 \\ 0.0000 &  0.0000 \end{array}\right] \la^2 
$$
$$
 G_c(\la)\approx
\left[\begin{array}{rr}  0.7641 & 0.3496 \\ 0.1274 & -0.5270 \end{array}\right] +
\left[\begin{array}{rr}  0.3496 &  0.4144 \\ -0.5270 &  0.6544 \end{array}\right] \la 
$$
and obtained a factorization residual norm
$ \|P(\la)-N_r(\la)G_c(\la) \|_F \approx  8.10^{-15}$.
The tolerance for the rank decisions was chosen as $tol = 1000 \epsilon_M$, where $\epsilon_M\approx 2.10^{-16}$ is the machine precision of the computer used.
To check the roots of $G_c(\la)$ we compared the (monic) characteristic polynomials of
$G_c(\la)$ and $G(\la)$ (the latter is $\la^2+\la-1$) and the error was approximately 
$2.10^{-15}$.

\begin{center}
    \begin{figure}
        \centering
        \includegraphics[scale=0.3]{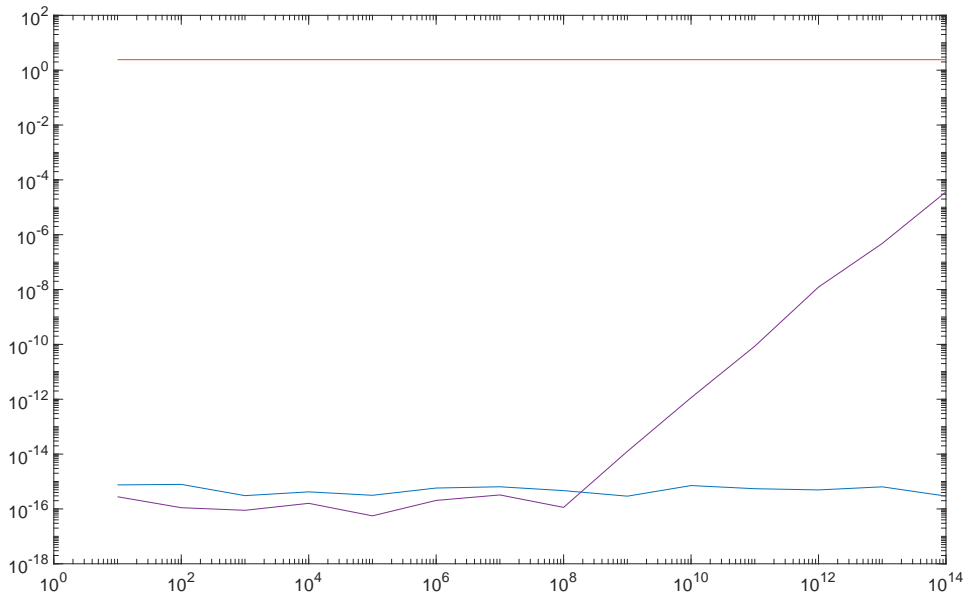}
        \caption{In blue, residual $\rho_1(k)$ for the computed GCRD factorization of \eqref{eq:secondexp} for varying values of $k$. In red, condition number $\rho_2(k)$ of the computed compact GCRD of \eqref{eq:secondexp}, evaluated at $\la=1$, for the same values of $k$. In purple, norm $\rho_4(k)$ of the vector $G_c'(0,k)e_1-G_c(0,k)e_2$ for the same values of $k$.}
        \label{fig:my_label}
    \end{figure}
\end{center}

As a second experiment, we constructed the input
\begin{equation}\label{eq:secondexp}
    P(\la,k)= Z \begin{bmatrix} \la^2 & 2 \la \\
0 & \la \\
\la & k \la + 1\\
0 & \la^2 \end{bmatrix} \in \C[\la]^{4 \times 2}, 
\end{equation} 
where $0 \leq k \in \R$ is a parameter, $Z \in \C^{4 \times 4}$ is a complex unitary matrix, randomly generated by first taking $X$ to be a realization of a $4 \times 4$ Gaussian random matrix, and then computing $Z$ as the unitary factor of the $QR$ factorization of $X$. It is not difficult to show that, for all $k, Z$ as above, any compact GCRDs of \eqref{eq:secondexp} must have the form $V(\la) \begin{bmatrix} \la & 1\\
0 & \la\end{bmatrix}$ where $V(\la) \in \C[\la]^{2 \times 2}$ is unimodular; in particular, all these GCRDs are regular, have a defective eigenvalue at $1$, and have a root polynomial \cite{DN,NV} at $0$ equal to $\begin{bmatrix} -1 & \la \end{bmatrix}^T$. On the other hand, when $k \rightarrow + \infty$, the columns of the input become unbalanced: the Frobenius norm of the second column diverges, whereas the Frobenius norm of the first column is independent of $k$ and equal to $\sqrt{2}$. Hence, when $k$ is large one could expect numerical difficulties. For $k=10,10^2,\dots,10^{14}$, we computed a GCRD factorization $P(\la,k)=N_r(\la,k) G_c(\la,k)$ and then four quantities: the factorization residual $\rho_1(k)=\|N_r(\la,k)G_c(\la,k)-P(\la,k)\|_F$, the quantity $\rho_2(k)=\kappa(G_c(1,k))$, where $\kappa$ is the $2$-norm condition number, and the norms $\rho_3(k)=\| G_c(0,k) e_1\|_2$ and $\rho_4(k)=\|G_c'(0,k) e_1 - G_c(0,k)e_2\|_2$. In particular, if $\rho_2(k)$ is not large then $G_c(\la,k)$ is correctly computed to be regular; and if $\rho_3(k)$ and $\rho_4(k)$ are both small then $\begin{bmatrix}-1 & \la \end{bmatrix}$ is indeed a root polynomial at $0$ for $G_c(\la)$. In our computations, $\rho_3(k)$ was computed to be exactly zero to machine precision for all values of $k$; Figure \ref{fig:my_label} below provides a logarithmic plot of $\rho_1(k)$, $\rho_2(k)$, and $\rho_4(k)$.
Note that, for $\log_{10}(k) \geq 9$, the root polynomials at the eigenvalue $0$ were not always accurately captured. Values of $k$ for which $\rho_4(k)$ is large suggest that the largest partial multiplicity of $0$ was computed to be $1$ in lieu of the correct value of $2$; but this behaviour in the presence of a Jordan chain of length $2$ is not surprising for an algorithm run in double precision.

The numerical tests discussed above are for small-sized and small-degree polynomial matrices. It should be noted that the problem of computing a GCRD is ill-posed. One way to see this is to observe, generically in the Euclidean topology, almost every set of input matrices only have right invertible (over $\R[\la]$) compact GCRDs. Recall \cite[Lemma 5.2]{nvdpen} that being right invertible over $\R[\la]$ is equivalent to having no finite eigenvalues and no left minimal indices. Thus, espcially when the size or the degree of the input grow larger, even a backward stable method can be expected to incorrectly approximate finite eigenvalues in the computed GCRD. Indeed, all the numerical experiments that we could find in the literature for (both exact and approximate) GCRD computation are for polynomial matrices of very small, say $O(1)$, size and degree. Nevertheless, in the next experiments we will test our code against higher values of the input size and of the input degree.

For our third experiment, we randomly generated ten examples of the form $P(\la)=M(\la)S(\la)N(\la) \in \R[\la]^{1000 \times 500}$ of normal rank $20$ where $M(\la)\in\R[\la]^{500\times 20}$ and $N(\la)\in\R[\la]^{20\times 500}$ are degree-$1$ polynomial matrices whose coefficients are realizations of Gaussian random matrices, and $S(\la)\in\R[\la]^{20\times 20} = I_{19} \oplus p(\la)$ where $p(\la)$ is a scalar polynomial of degree $4$ with coefficients that are realizations of normal random variables. Then, almost surely, the matrices $P(\la)\in\R[\la]^{1000\times 500}$
are of degree $6$ and have as finite Smith zeros the roots of the
fourth degree polynomial $p(\la)$. Moreover, again with probability $1$, by Remark \ref{rem1} any compact GCRD $G_c(\la) \in \R[\la]^{20 \times 500}$ has $4$ finite eigenvalues and $480$ right minimal indices. For these larger size inputs, we set the tolerance for rank decisions at $tol=10^4 \epsilon_M$. Table \ref{table1} below yields, for these randomly generated examples, the norms of the computed factors $N_r(\la)$ and $G_c(\la)$ of $P(\la)$, the norm of the residual $Res(\la):=P(\la)-N_r(\la)G_c(\la)$ and the inverse condition numbers $\kappa(G_c(\la_i))^{-1}:=
\sigma_r(G_c(\la_i))/\sigma_1(G_c(\la_i))$ of the 
matrices $G_c(\la_i)$ evaluated in the Smith zeros of $\la_i$ of $P(\la)$. 

\medskip

\noindent

\begin{center}
\begin{table}[h!]
    \centering
   \begin{tabular}{rrrrrrr}
$\|N_r(:)\|_2$ & $\|G_c(:)\|_2$ & $\|Res(:)\|_2$ & $\kappa(G_c(\la_1))^{-1}$ &  $\kappa(G_c(\la_2))^{-1}$ &  $\kappa(G_c(\la_3))^{-1}$ &  $\kappa(G_c(\la_4))^{-1}$ \\ \hline
 1.0013e+00   &4.4721e+00 &  2.9624e-15 &  1.7509e-15  & 1.7509e-15  & 2.7147e-15  & 2.7147e-15\\
   1.0003e+00   &4.4721e+00  & 3.7683e-15&   1.6676e-16  & 3.3772e-16&   3.3772e-16   &4.1002e-15\\
   9.9912e-01   &4.4721e+00 &  3.2400e-15 &  3.0898e-16 &  1.0641e-15 &  4.1927e-15   &4.1927e-15\\
   9.9966e-01   &4.4721e+00   &4.4029e-15  & 2.5095e-16&   1.9051e-15  & 1.3357e-15   &1.3357e-15\\
   9.9992e-01   &4.4721e+00  & 3.0857e-15   &3.3193e-16   &9.9503e-16  & 9.9503e-16   &2.1699e-15\\
   1.0003e+00   &4.4721e+00 &  3.8058e-15 &  7.0299e-16  & 1.1784e-15  & 1.1784e-15   &5.3952e-15\\
   1.0007e+00   &4.4721e+00   &4.6801e-15  & 2.1369e-16   &2.1369e-16  & 6.9481e-16   &6.9481e-16\\
   1.0005e+00   &4.4721e+00   &4.0644e-15   &1.4874e-16  & 8.0534e-16  & 9.3845e-16   &9.3845e-16\\
   9.9994e-01   &4.4721e+00  & 3.3280e-15&   3.0311e-16 &  3.0311e-16  & 7.7456e-16   &2.4649e-15\\
   1.0005e+00  & 4.4721e+00 &  6.4201e-15 &  6.5535e-16&   6.4943e-16  & 4.9751e-15   &7.6166e-15\\
\end{tabular}
    \caption{Norms of the factors $N_r(\la)$, $G_c(\la)$ and the residual matrix $R(\la)$, and condition numbers of $G_c(\la)$ at the eigenvalues $\la_i,i=1,\ldots,4$, for ten randomly generated compound matrices $P(\la)$ as in \eqref{eq:compound} with $m=1000$ rows, $n=500$ columns, normal rank $r=20$, degree $6$, and $4$ finite eigenvalues.}
    \label{table1}
\end{table}

\end{center}

It can be seen from Table \ref{table1} that in all of the 10 examples the algorithm found the correct normal rank $r=20$ and a right factor $G_c(\la)$ that contains the four Smith zeros of $P(\la)$; the accuracy in the computation of the eigenvalues is remarkable and even somewhat surprising for an input of $O(10^3)$ size, given the remarks above on the ill-posedness of the problem. Moreover, the norms of the two factors $N_r(\la)$ and $G_c(\la)$ are very reasonable, which can be explained from the use of orthogonal transformations
applied to the system matrix pencil. As pointed out in Remark \ref{rem:boundG} the Frobenius norm of 
$G_c(\la)$ is approximately the square root of its number of rows; note that $\sqrt{20} \simeq 4.4721$, so the algorithm has always computed the normal rank correctly. We conclude that, in the case of small input degree but relatively large input size, the experimental results are thus compatible with a backward stable algorithm, both for the compact GCRD extraction and for the computation of the factorization $P(\la)=N_r(\la)G_c(\la)$.

\medskip
For our fourth experiment, we repeated a similar numerical test as the previous one but with different parameters, to test the behaviour of algorithm for moderate-degree input. This time $P(\la)=M(\la)S(\la)N(\la) \in \R[\la]^{4 \times 3}$ is constructed by randomly generating degree-$10$ $M(\la)\in \R[\la]^{4 \times 2}$ and $N(\la)\in \R[\la]^{2 \times 3}$, and $S(\la)=1\oplus p(\la)$ where $p(\la)$ is a randomly generated scalar polynomial of degree $4$. Hence, almost surely, $P(\la) \in \R[\la]^{4 \times 3}$ has degree $24$ and normal rank $2$. Therefore, almost surely, any compact GCRD $\R[\la]^{2 \times 3}$ has $4$ eigenvalues and $1$ right minimal index. The same residual tests as before have been computed and are reported in Table \ref{table2} below.
\medskip
\begin{center}
\begin{table}[h!]
    \centering
    \begin{tabular}{rrrrrrr}
$\|N_r(:)\|_2$ & $\|G_c(:)\|_2$ & $\|Res(:)\|_2$ & $\kappa(G_c(\la_1))^{-1}$ &  $\kappa(G_c(\la_2))^{-1}$ &  $\kappa(G_c(\la_3))^{-1}$ &  $\kappa(G_c(\la_4))^{-1}$ \\ \hline
1.0300e+00 &  1.4142e+00&   1.3774e-13&   7.6420e-16&   7.6420e-16&   4.7405e-16&   3.9246e-13\\
   1.0052e+00&   1.4142e+00&   1.4404e-13&   2.9781e-13&   5.4190e-13  & 3.4446e-13&   3.4446e-13\\
   5.2905e+00&   1.7321e+00 &  7.2512e-05   &2.4736e-15 &  3.8403e-16  & 3.8403e-16&   1.5077e-14\\
   1.6320e+00&   1.7321e+00  & 3.1000e-13  & 1.3971e-16&   4.3847e-15  & 4.3847e-15&   5.2610e-14\\
   4.2221e+00&   1.7321e+00   &1.2369e-08 &  6.9758e-15&   6.9758e-15  & 2.3524e-15&   1.4909e-14\\
   6.7245e+00&   1.7321e+00 &  9.1725e-05&   6.1401e-16 &  6.1401e-16  & 2.3190e-14&   2.3190e-14\\
   9.5240e-01&   1.4142e+00&   2.7136e-13   &5.5111e-16  & 1.6181e-15  & 1.6181e-15&   1.0418e-12\\
   1.0000e+00&   1.7321e+00 &  2.9692e-13  & 1.0000e+00   &1.0000e+00  & 1.0000e+00&   1.0000e+00\\
   1.0000e+00&   1.7321e+00  & 2.4253e-14 &  1.0000e+00&   1.0000e+00 &  1.0000e+00&   1.0000e+00\\
   1.0295e+00&   1.4142e+00   &8.7189e-14&   3.8456e-15 &  3.9381e-14&   3.9381e-14 &  1.3311e-07\\
\end{tabular}
    \caption{Norms of the factors $N_r(\la)$, $G_c(\la)$ and the residual matrix $R(\la)$, and condition numbers of $G_c(\la)$ at the eigenvalues $\la_i,i=1,\ldots,4$, for ten randomly generated compound matrices $P(\la)$ as in \eqref{eq:compound} with $m=4$ rows, $n=3$ columns, normal rank $r=2$, degree $24$, and $4$ finite eigenvalues.}
    \label{table2}
\end{table}
\end{center}
\medskip
Table \ref{table2} shows that, when the input degree grows, computing a compact GCRD correctly becomes more challenging. For example, the factorization residual is not always small, indicating that, when the input degree is $\gg 1$, then the algorithm is not backward stable for the computation of the factorization $P(\la)=N_r(\la)G_c(\la)$. The results are, however, still compatible with the algorithm being backward stable as a method to compute $G_c(\la)$ only, even though the analysis to explain this is delicate. First, recalling Remark \ref{rem:boundG} and noting that $\sqrt{2} \simeq 1.4142$ and $\sqrt{3} \simeq 1.7321$, it is clear that most of the times the algorithm has overestimated the normal rank of the input; occasionally, this has also led to a high factorization residual. It also becomes relatively common to inaccurately capture the finite eigenvalues that the GCRD should (in exact arithmetic) have, or even to miss them entirely.  As mentioned above, these shortcomings are not an indication of backward instability of the method as an algorithm to compute a compact GCRD. Indeed, it is known \cite{dd} that polynomial matrices with full column rank and no eigenvalues are everywhere dense (in the Euclidean topology) in $\R[\la]_{24}^{4 \times 3}$, the vector space of $4 \times 3$ real polynomial matrices of degree at most $24$. Hence, most small perturbation of the input are such that their compact GCRD is a $3 \times 3$ unimodular matrix. In fact, the lines of the table above where the inverse condition number at the expected eigenvalues is $1$ are caused by the fact that the algorithm computed $G_c(\la)$ to be a constant orthogonal matrix: by the discussion above, this is still compatible with backward stability, because there exist (infinitely many) polynomial matrices whose GCRD is a constant unitary matrix and that are arbitrarily close to the actual input. It is in some sense remarkable that, at least in some cases, the correct eigenvalues are not missed even if the degree of the input is $\gg 1$. 

Finally, we tested the effect of the tolerance with which the algorithm runs (which is a crucial parameter in particular for the staircase algorithm): in Table \ref{table3}, we report the same tests for the \emph{same} inputs (row by row) as Table \ref{table2}, but having re-run the algorithm with a tolerance increased by a factor $10$.

\begin{center}
\begin{table}[h!]
    \centering
    \begin{tabular}{rrrrrrr}
$\|N_r(:)\|_2$ & $\|G_c(:)\|_2$ & $\|Res(:)\|_2$ & $\kappa(G_c(\la_1))^{-1}$ &  $\kappa(G_c(\la_2))^{-1}$ &  $\kappa(G_c(\la_3))^{-1}$ &  $\kappa(G_c(\la_4))^{-1}$ \\ \hline
1.0300e+00 &  1.4142e+00&   1.3774e-13 &  7.6420e-16 &  7.6420e-16 &  4.7405e-16  & 3.9246e-13\\
   1.0052e+00 &  1.4142e+00 &  1.4404e-13&   2.9781e-13&   5.4190e-13&   3.4446e-13  & 3.4446e-13\\
    1.1276e+00  & {\bf 1.4142e+00}&   {\bf 4.9997e-12} &  1.2706e-13   & 6.3114e-13 &   6.3114e-13  & 4.2452e-12\\
   9.1563e-01  & 1.7321e+00   &{\it 5.2990e-12}  & 5.8681e-16  & 3.0430e-15  & {\bf 3.0430e-15}  & 1.5136e-13\\
   4.2221e+00  & 1.7321e+00  & 1.2369e-08   &6.9758e-15 &  6.9758e-15  & 2.3524e-15  & 1.4909e-14\\
   1.0362e+00  & {\bf 1.4142e+00} &  {\bf 7.6057e-13}&   1.4496e-14&   1.4496e-14  & 1.4221e-12  & 1.4221e-12\\
   9.5240e-01  & 1.4142e+00   &2.7136e-13 &  5.5111e-16   &1.6181e-15  & 1.6181e-15  & 1.0418e-12\\
   7.6350e+01  & 1.7321e+00  & {\it 3.8364e-06}  & {\bf 6.4356e-17}  & {\bf 1.7544e-14}  & {\bf 5.7841e-14}  & {\bf 5.7841e-14}\\
   9.9824e-01  & {\bf 1.4142e+00} &  {\it 5.7547e-12}&   {\bf 6.7049e-15} &  {\bf 1.9289e-11}  & {\bf 3.5895e-11}  & {\bf 3.5895e-11}\\
   1.0295e+00  & 1.4142e+00&   8.7189e-14 &  3.8456e-15&   3.9381e-14  & 3.9381e-14  & 1.3311e-07\\
\end{tabular}
    \caption{Norms of the factors $N_r(\la)$, $G_c(\la)$ and the residual matrix $R(\la)$, and condition numbers of $G_c(\la)$ at the eigenvalues $\la_i,i=1,\ldots,4$, for the same inputs as Table \ref{table2}, but having run the algorithm with a tolerance $10$ times higher. Results that have improved significantly after the increase in tolerance are in bold font; those that have worsened significantly are in italic.}
    \label{table3}
\end{table}
\end{center}
It can be seen in Table \ref{table3} that the numerical estimations of the normal rank and of the finite eigenvalues often improve by increasing the tolerance; however, occasionally this comes at the expense of a higher factorization residual. How to set the tolerance a priori with respect to the known parameters (such as size and degree of an input) is generally a difficult problem and left for future research.
\section{Concluding remarks}  \label{sec:conclusions}

In this paper we revisited the problem of Greatest Common Right Divisor extraction in its most general form, where the given matrix $P(\la)$ does not need to have full column rank, and where the GCRD need not be square. Along the way, we introduced the so-called compact GCRD's which play a central role here since all other GCRDs can be derived from them. Moreover, we presented a novel algorithm that computes such a minimal size GCRD. The compact factorization is expected to play an important role in other problems such as the computation of a compact Smith-like decomposition $P(\la)= U_r(\la) S_r(\la) V_r(\la)$, where the middle factor $S_r(\la)$ is regular and is of dimension $r\times r$.

\end{document}